\documentclass[11pt]{article}
\usepackage{style}
\usepackage{booktabs}
\newcommand{\ManuscriptTableFont}{\small}
\bibpunct[, ]{(}{)}{,}{a}{}{,}
\newcommand{\R}{\mathbb{R}}
\newcommand{\E}{\mathbb{E}}
\newcommand{\Prob}{\mathbb{P}}
\newcommand{\TV}{\mathrm{TV}}
\newcommand{\Var}{\mathrm{Var}}

\newcommand{\ind}{\mathbf{1}}

\newcommand{\nulla}{\varnothing}

\title{\LARGE\bfseries Gradient Estimation for LLM-Embedded Simulation Optimization}
\author{Ruiyang Jin, Siyang Gao\\[0.4em]
{\small \textit{Department of Systems Engineering, City University of Hong Kong, Hong Kong, China}}}
\date{\vspace{-0.4in}}
\hypersetup{pdftitle={Gradient Estimation for LLM-Embedded Simulation Optimization},pdfauthor={Ruiyang Jin, Siyang Gao}}

\begin{document}
\maketitle

\begin{abstract}
Large language models (LLMs) are increasingly embedded in operational simulators to represent context-dependent human decisions, which opens new possibilities for simulation optimization. However, the LLM acts as a black box, so derivatives through its decisions are unavailable for gradient-based policy search. Finite difference estimation can avoid this requirement, but the required paired simulations that sample LLM actions independently often disagree, and such action switches persist under small perturbations, which inflates the gradient variance to order $c^{-2}$ for perturbation radius $c$. To control this variance, we develop a behavioral common random number framework. It first applies maximal coupling to the LLM's action distributions, preserving each simulation's behavioral law while making action switches as rare as possible, which reduces the variance bound to $O(c^{-1})$. It then applies first-switch conditioning, which averages over whether a switch occurs, together with an impact-aware extension that prioritizes influential switches. Under suitable stability conditions, both keep the variance bounded as $c$ shrinks. We further establish improved complexity bounds for projected stochastic approximation. Experiments on synthesized and LLM-embedded simulators confirm the theory and show substantial performance gains.
\end{abstract}


\section{Introduction}
\label{sec:introduction}

Simulation optimization is used in operations research to improve policies whose performance has no closed-form expression and must be estimated through stochastic simulations \citep{fu2015handbook}. In many operational simulators, system evolution depends on both exogenous inputs, such as arrivals and service times, and decisions made by agents within the simulated system \citep{nelson2021foundations}. In service operations, for example, call-center performance depends on customers' patience, abandonment, and responses to waiting-time or callback information \citep{gans2003telephone,yu2023managing}. In inventory systems, stockout performance depends on whether customers substitute, backorder, or cancel after observing product and fulfillment information \citep{honhon2010assortment,chen2020dynamic}. These behaviors have traditionally been characterized by parametric choice models, hazard models, or rule-based response modules. Such models require specifying the state variables and explicit behavioral kernels in advance \citep{ye2020hazard}. Moreover, such models generally cannot process multimodal contexts, such as natural language and images.

Large language models (LLMs) offer another way to represent context-dependent decisions. They can learn in context from textual prompts \citep{brown2020language}, perform multi-step reasoning with appropriate prompting \citep{wei2022chain}, and generate actions in interactive environments \citep{yao2023react}. These capabilities have led to their use as behavioral surrogates in computational social, economic, and user-interaction systems \citep{park2023generative,goli2024frontiers}. A natural next step is to embed LLMs in operational simulators. The simulator retains its conventional physical dynamics, while the embedded LLM makes behavioral decisions from the contexts. This emerging practice is highly relevant and critically important to simulation optimization, because policies can then be evaluated and improved against behaviors that respond to rich contexts, without specifying a behavioral model in advance.

Embedding an LLM, however, also changes how the simulator generates its outputs. The sampled decisions affect later dynamics and create feedback between behavioral and physical evolution. Moreover, the LLM typically acts as a black box. Under a sampled-action interface, the simulator submits a context and observes an action, while the probabilities, internal randomization, and derivatives may be unavailable.

Gradient estimation supports policy search in these simulators, but the behavioral module creates difficulties for direct gradient methods. Pathwise and likelihood-ratio methods require sample-path regularity or information about the behavioral probabilities and their derivatives \citep{glasserman1991gradient,nelson2021foundations,mohamed2020monte}. Finite-difference estimators remain applicable because they use simulation outputs \citep{larson2019derivative,berahas2022theoretical}. Common random numbers (CRNs) can reduce their variance by sharing physical primitives to couple evaluations at nearby policies \citep{glasserman1992guidelines}. However, this leaves the behavioral calls independent, so the two evaluations can draw different actions at the same event. Such action switches do not vanish as the perturbation shrinks, and the resulting gradient estimates can remain unstable despite shared physical randomness.

We study this problem as \emph{LLM-embedded simulation optimization} and develop a behavioral common random number (BCRN) construction for coordinate-wise finite-difference gradient estimation (CGE). Just as CRNs share physical randomness, BCRN couples the behavioral draws of nearby policy evaluations while preserving their marginal laws. Implementing this coupling requires access beyond sampled actions. Our score-access interface evaluates each context separately and converts the returned action scores into explicit probability vectors. Maximal coupling then samples the paired actions from these two probability vectors so that they agree as often as possible. This preserves the marginals exactly and requires no derivatives of the LLM or downstream simulator.

Maximal coupling reduces the probability of behavioral disagreements, but a rare disagreement can still substantially change the trajectory. Because the finite difference divides this change by a small perturbation radius, rare switches remain a major source of noise. We apply first-switch conditioning to address this remaining noise. Rather than leaving it to chance whether a replication contains a switch, it averages over the no-switch and first-switch cases. Its impact-aware extension further uses an impact predictor to select the first-switch event, which permits allocation toward influential branches and further reduces the variance.

Our contributions are as follows.
\begin{itemize}
\item We identify how behavioral randomness destabilizes finite-difference gradient estimates in LLM-embedded simulation optimization even when all physical randomness is shared. Under certain conditions, CGE with independent sampling (I-CGE) has a variance floor of order $c^{-2}$, where $c$ is the perturbation radius. This floor arises because paired behavioral actions can disagree with non-vanishing probability as \(c\) decreases.
\item We introduce CGE with maximal-coupling sampling (MS-CGE) to implement BCRN through a score-access interface. It makes action distributions explicit while treating the LLM and downstream simulator as black boxes. Maximal coupling preserves behavioral marginals exactly and minimizes the probability of action switches. Under locally Lipschitz pre-switch stability, MS-CGE obtains an improved variance bound of $O(c^{-1})$.
\item To control noise from rare but potentially influential action switches, we propose CGE with first-switch conditioning (FS-CGE) and the impact-aware variant (IA-FS-CGE). Both estimators preserve the central-difference mean, and IA-FS-CGE can further reduce variance through impact-aware selection. Under pathwise Lipschitz stability, their variance is uniformly bounded in \(c\). Using these variance bounds, we establish improved complexity bounds for applying the proposed estimators in projected stochastic approximation.
\item We evaluate the variance-reduction mechanism and optimization performance under common budgets. Synthesized experiments illustrate the variance behavior of the four gradient estimators. Experiments on a software-onboarding problem with an LLM-embedded simulator show the same variance trends and validate the proposed estimators' optimization performance.
\end{itemize}

\paragraph{Related work.}
Our work draws on research in LLM-based behavioral simulation, which provides a flexible way to represent decisions driven by agent characteristics and situational context, including interactions that generate emergent collective behavior \citep{park2023generative}. Applications span survey responses and experimental choices \citep{argyle2023out}, economic decision-making \citep{goli2024frontiers}, and user feedback for recommender and conversational systems \citep{yoon2024evaluating, zhang2025llmuser}. The common emphasis is on constructing context-dependent agents and evaluating their behavior within simulated systems \citep{karetnikov2026proxies, ghaffarzadegan2024generative}. Existing literature shows LLMs' potential as context-dependent behavioral agents, but focuses primarily on behavioral realism, population alignment, and system-level evaluation.

Stochastic approximation searches over continuous policies by updating them with noisy gradient estimates, following the Robbins-Monro and Kiefer-Wolfowitz procedures \citep{robbins1951stochastic,kiefer1952stochastic}. Its efficiency depends on gradient accuracy and evaluation cost. The available estimators include direct methods, such as pathwise and likelihood-ratio estimation, and finite differences \citep{fu2015handbook}. Pathwise methods, including infinitesimal perturbation analysis (IPA), differentiate a realized trajectory under sample-path regularity \citep{glasserman1991gradient,LEcuyer1990}. However, discrete actions can violate this regularity. With the primitive randomness fixed, an action may remain unchanged under small policy perturbations and then jump at a decision boundary. As a result, its pathwise derivative can miss behavioral sensitivity or be undefined at a boundary, even though the expected objective changes with the action probabilities. Likelihood-ratio methods instead require derivatives of the log probabilities of sampled actions \citep{glynn1990likelihood,fu2006gradient}. Our interface supplies probability values for coupling, but neither these derivatives nor derivatives through the simulator.

Finite-difference gradient estimation approximates derivatives by comparing objective values at perturbed policies \citep{berahas2022theoretical}. CGE applies a two-sided perturbation to each coordinate and assembles the resulting partial-derivative estimates \citep{kiefer1952stochastic}. CRNs can reduce their noise by increasing covariance between the paired outputs \citep{glasserman1992guidelines}. In classical discrete-event simulation, random streams are synchronized according to their operational roles \citep{asmussen2007stochastic}. In an LLM-embedded simulator, physical CRNs can align exogenous inputs but do not by themselves align categorical behavioral sampling.

Coupled finite differences provide such variance reduction when local transition probabilities can be evaluated but their derivatives are unavailable. In structured Markov models, suitable couplings give $O(c^{-1})$ variance under regularity conditions \citep{anderson2012efficient}. For categorical transitions, Gumbel-max couples draws from known probability vectors without requiring an analytic mapping from contexts to probabilities \citep{oberst2019counterfactual}. Learned coupling mechanisms can also reduce the variance of downstream contrasts \citep{lorberbom2021learning}. These results motivate using probability values to coordinate behavioral sampling while leaving the context dynamics undifferentiated.

Conditional Monte Carlo averages over discrete branches when their probabilities are tractable and conditional continuations can be simulated \citep{fu1997conditional}. Related methods include measure-valued differentiation for Markov chains \citep{heidergott2008measure} and stochastic derivatives that propagate weighted alternative executions through discrete randomness \citep{arya2022automatic}. First-divergence decompositions appear in phantom estimators and GradInf \citep{fu1997conditional,arya2026gradinf}. GradInf also combines residual inference with stratified resampling or twisted sequential Monte Carlo to favor current or predicted downstream perturbations. Second-moment sampling distributions and regression-then-sampling designs are established for stochastic simulators \citep{chen2019importance}.

The remainder of the paper follows this organization. Section~\ref{sec:problem_formulation} introduces the finite-horizon simulator and behavioral kernel. Section~\ref{sec:gradient_estimation} presents I-CGE, the score--softmax construction for MS-CGE, FS-CGE and its impact-aware extension. Section~\ref{sec:analysis_PGE} establishes their bias, variance, and optimization-complexity results. Sections~\ref{sec:simulation} and~\ref{sec:conclusion} report the numerical experiments and conclusions.

\section{Problem Formulation}\label{sec:problem_formulation}

\noindent\textit{Notation.} We use $\E$, $\Prob$, and $\Var$ for expectation, probability, and variance, and $\ind\{\cdot\}$ for the indicator of an event. For vectors, $\langle\cdot,\cdot\rangle$ and $\|\cdot\|$ denote the Euclidean inner product and norm, while $\|\cdot\|_\infty$ is the maximum absolute component. The vector $e_j$ is the $j$th standard basis vector in $\R^d$, and $\partial_j J$ denotes the $j$th component of the gradient $\nabla J$. $\Pi_\Theta$ denotes the Euclidean projection onto the set $\Theta$.

We consider a simulator with finite horizon $T<\infty$ and $m<\infty$ potential behavioral events. Let $\Theta\subseteq\R^d$ be a closed convex set of policy centers, and let $J$ be the expected objective evaluated by simulation. At selected events, an LLM-based behavioral module generates a finite-valued action that affects the subsequent physical evolution. For any $t\leq T$, let $X_t^\theta$ denote the physical state at time $t$ under policy $\theta$. The state evolves as
\begin{equation*}
X_{t+1}^{\theta}=\mathcal{T}\big(X_t^{\theta},\theta,\xi_t, \omega_t \big),
\qquad t=0,1,\dots,T-1,
\end{equation*}
where $\mathcal{T}$ is the physical state transition operator, $\xi_t$ is exogenous physical randomness, and $\omega_t$ records the realized behavioral actions used at time $t$.

We label event $e\in\{1,\dots,m\}$ with a simulator-defined identifier, for example, a tuple of the form
$\text{(entity id, decision type, occurrence number)},$
and use $t_e^{\theta}$ to denote the time index of event $e$.
At event $e$, the simulator forms a context $C_e^{\theta}=\Gamma_e\big(H_e^\theta,\zeta_e\big)\in\mathcal C$, where $H_e^\theta$ is the history up to event $e$ and $\zeta_e$ contains persistent entity attributes, such as customer history, a patient profile, or order characteristics. The context may combine multimodal information.

At each event $e=1,\ldots,m$, the LLM-based module uses a context-dependent behavioral kernel $P_{b,e}$ on the finite action set $\mathcal A$ to draw $A_e^\theta\sim P_{b,e}(\cdot\mid C_e^\theta)$.
The prompt, LLM, decoder, and action parser together may determine this kernel, which can vary across behavioral events. Under sampled-action access, the simulator observes the realized action but cannot inspect its probability vector or sampling seed. Therefore, coordinating the behavioral draws of two trajectories requires the additional score access introduced in Section~\ref{subsec:maximal-fd-estimator}.

The goal of simulation optimization is to minimize the expected loss
\begin{equation*}
J(\theta)=\E[L(\theta, \boldsymbol\xi, \boldsymbol\omega)],
\end{equation*}
where $L(\theta, \boldsymbol\xi, \boldsymbol\omega)$ is the integrable total trajectory loss. The random paths over the trajectory are $\boldsymbol\xi=(\xi_0,\ldots,\xi_{T-1})$ and $\boldsymbol\omega=(\omega_0,\ldots,\omega_{T-1})$. The physical primitives include the initial state and entity attributes, so the expectation also averages over these inputs. 

\begin{example}
\label{ex:inventory}
Consider a retailer that manages $d$ products over $T$ periods with continuous base-stock levels $\theta\in\Theta\subseteq\R^d$. The physical state $X_t^\theta$ records on-hand inventory, outstanding orders, and backorders, and $\xi_t$ contains customer arrivals, requested products, and replenishment lead times. When a requested product is out of stock, a behavioral event occurs at the customer's arrival period $t_e^\theta$. Since customers are generated in advance, the event of customer $k$ is labeled by (customer $k$, stockout response, 1) and $m$ is finite. How a customer responds depends on brand loyalty, urgency, the perceived similarity of the substitute, and the wording of the delivery promise, which are difficult to capture with a parametric choice model with fixed substitution rates. We therefore let an LLM act as the customer. The context $C_e^\theta$ is a prompt that combines the customer's profile $\zeta_e$, such as purchase history and preferences, with a product page showing the requested item, a substitute, and the promised delivery date. The LLM's sampled reply is parsed into an action in $\mathcal A=\{\text{substitute},\text{backorder},\text{cancel}\}$. A substitution depletes the substitute's inventory, and a backorder consumes future replenishment, so each action changes the prompts faced by later customers. The loss $L$ charges holding, backorder, and lost-sales costs. Thus $\theta$ affects $J$ directly through holding costs and indirectly through the prompts, while the LLM's response probabilities can only be sampled and cannot be differentiated with respect to $\theta$. If a product is in stock under one policy but not under another, the event is active only under the latter and receives the null action under the former, as described in the remark below.
\end{example}

\begin{assumption}\label{assump:smoothness}
There exist constants $L_J>0$ and $c_0>0$ such that the simulator and $J$ are defined on the closed $c_0$-neighborhood of $\Theta$ and, for every two points $\theta,\theta'$ in this neighborhood,
\[
\left\|\nabla J(\theta')-\nabla J(\theta)\right\|\le L_J \|\theta'-\theta\|.
\]
\end{assumption}
Assumption~\ref{assump:smoothness} is the standard Lipschitz-gradient condition. The neighborhood condition permits two-sided perturbations without clipping.
\begin{remark}
The paired estimators require a common causal ordering of at most $m$ potential behavioral events. To align events across trajectories, we assign a deterministic null action $\nulla$ whenever an event is inactive. This convention preserves event-label alignment when an event is active under one policy but inactive under the other. Without loss of generality, the same candidate action set $\mathcal A$ applies to all events after including the null action. Moreover, physical primitives are independent of the primitive behavioral sampling seeds and may be generated in advance.
\end{remark}

\section{Gradient Estimation}\label{sec:gradient_estimation}
We use two-sided CGE because discrete actions in LLM-embedded simulations can violate IPA regularity \citep{glasserman1991gradient,arya2022automatic} and the probability derivatives needed by likelihood-ratio methods are unavailable \citep{glynn1990likelihood,mohamed2020monte}. The estimators use simulation outputs without differentiating the behavioral module or simulator.

For a policy $\theta\in\Theta$, coordinate $j\in\{1,2,\cdots, d\}$, and radius $0<c\le c_0$, write $\theta_j^+:=\theta+ce_j$ and $\theta_j^-:=\theta-ce_j$. The central difference
\[
D_{j,c} J(\theta)=\frac{J(\theta_j^+)-J(\theta_j^-)}{2c}
\]
approximates $\partial_j J(\theta)$. The following methods estimate $D_{j,c}J(\theta)$ by evaluating $J(\theta_j^\pm)$ with different ways of sampling the paired trajectories. Superscripts $\mathrm{I}$, $\mathrm{MS}$, $\mathrm{FS}$ and $\mathrm{IA}$ denote the four estimators.

\subsection{CGE with Independent Sampling}\label{subsec:independent-fd}
I-CGE samples behavioral actions independently in the two perturbed systems. Let index $i$ denote an independent paired replication. For coordinate $j$ and replication $i$, the systems share the physical primitives $\boldsymbol\xi_{j,i}$ but query the behavioral module independently at their evolving contexts. Their joint action distribution at event $e$ is
\begin{equation*}
Q_{e}^{\mathrm{I}}(\cdot,\cdot\mid C_{e,j,i}^+,C_{e,j,i}^-)
=P_{b,e}(\cdot\mid C_{e,j,i}^+)\otimes
P_{b,e}(\cdot\mid C_{e,j,i}^-).
\end{equation*}
This product describes two independent black-box calls. The system need not observe their probability vectors or control their internal randomization.

Let $L_{j,i}^{\mathrm{I},\pm}$ be the trajectory losses, so $\frac{L_{j,i}^{\mathrm{I},+}-L_{j,i}^{\mathrm{I},-}}{2c}$ is one contribution to the $j$th partial gradient estimate. Averaging $n$ replications gives the estimate for coordinate $j$ and the full gradient estimate:
\begin{equation*}
\widehat G_{j,n}^{\mathrm{I}}(\theta,c)
:=\frac{1}{n}\sum_{i=1}^n
\frac{L_{j,i}^{\mathrm{I},+}-L_{j,i}^{\mathrm{I},-}}{2c},\qquad
\widehat{\boldsymbol g}_{n}^{\mathrm{I}}(\theta,c)
:=\sum_{j=1}^d\widehat G_{j,n}^{\mathrm{I}}(\theta,c)e_j.
\end{equation*}
As a result, a full I-CGE gradient uses $dn$ trajectory pairs. Each side follows its standalone behavioral kernel, giving each component mean $D_{j,c}J(\theta)$.

Independent behavioral calls can disagree with non-vanishing probability as the policy perturbation shrinks. Under the uniform nondegeneracy and loss-separation conditions, $\Var(\widehat G_{j,n}^{\mathrm I}(\theta,c))\ge v_0/(nc^2)$ for a constant $v_0>0$ (see Theorem~\ref{thm:indfloor} and Corollary~\ref{cor:indmse}). We provide the formal conditions and proof in Appendix~\ref{app:proof-indfloor}. The floor persists despite shared physical randomness and implies $\Omega(n^{-1/2})$ MSE under the radius schedule $c_n=\Theta(n^{-1/4})$.

\subsection{CGE with Maximal-Coupling Sampling}\label{subsec:maximal-fd-estimator}
To reduce action disagreement from independent sampling, CGE with maximal-coupling sampling (MS-CGE) couples behavioral actions between the perturbed systems. Since sampled-action access does not expose the required probabilities, MS-CGE uses a score-access interface: it constructs a distribution from the scores at each context, then samples a joint law with these distributions.

We say two contexts are \emph{matched} when they correspond to the same behavioral event, simulated entity, and immutable entity attributes.
For each context $C\in\mathcal C$, the LLM returns the score vector $s_e(C)=\bigl(s_{e,a}(C):a\in\mathcal A\bigr)\in\R^{|\mathcal A|}$ instead of a sampled action. Its component $s_{e,a}(C)$ measures the relative preference for action $a$ in that context. Given a fixed temperature $\tau>0$, MS-CGE converts the scores into the marginal behavioral kernel
\begin{equation}
P_{b,e}(a\mid C) =
\frac{\exp\{s_{e,a}(C)/\tau\}}
{\sum_{b\in\mathcal A}\exp\{s_{e,b}(C)/\tau\}},
\qquad a\in\mathcal A,
\label{eq:score-softmax-kernel}
\end{equation}
where temperature $\tau$ controls concentration. Smaller values of $\tau$ amplify score differences, and the denominator normalizes the probabilities. We fix the LLM, its prompt or feature construction, and $\tau$ during gradient estimation and optimization. All methods below use this constructed law $P_{b,e}$, so marginal preservation refers to the same behavioral model throughout the comparison. 

We define a \textit{pair oracle} as a joint action kernel whose two marginals equal the prescribed context-specific distributions. On a finite action space, the minimum disagreement probability over all such couplings is $\TV(\mu,\nu)=\tfrac12\sum_a|\mu(a)-\nu(a)|$, and a maximal coupling attains this minimum \citep[Proposition~4.7]{levin2017markov}. For matched contexts $(C_e^+,C_e^-)$, define
\[
p_e^\pm(a):=P_{b,e}(a\mid C_e^\pm),\qquad
h_e(a):=\min\{p_e^+(a),p_e^-(a)\},\qquad
\rho_e:=\sum_{a\in\mathcal A}h_e(a),
\]
where the overlap $h_e(a)$ is the mass that both systems assign to action $a$, and $\rho_e$ is the total mass available for a common draw.
The MS-CGE pair oracle is defined as follows.

\begin{definition}
\label{def:external-maximal}
If $\rho_e<1$, the external pair oracle uses the joint law
\begin{equation*}
Q_e^{\mathrm{MS}}(a,b\mid C_e^+,C_e^-)
=h_e(a)\ind\{a=b\}
+\frac{(p_e^+(a)-h_e(a))(p_e^-(b)-h_e(b))}{1-\rho_e}.
\end{equation*}
If $\rho_e=1$, it draws $A_e\sim p_e^+=p_e^-$ and sets $A_e^+=A_e^-=A_e$. 
\end{definition}

When $\rho_e<1$, the sampler draws from the diagonal overlap with probability $\rho_e$ and otherwise draws from the two normalized residual distributions.
The residual distributions contain each side's remaining mass and have disjoint supports. Writing $Q_{e,+}^{\mathrm{MS}}$ and $Q_{e,-}^{\mathrm{MS}}$ for the marginals of $Q_e^{\mathrm{MS}}$ gives
\begin{equation}
Q_{e,+}^{\mathrm{MS}}=p_e^+,\qquad
Q_{e,-}^{\mathrm{MS}}=p_e^-,\qquad
\Prob(A_e^+\ne A_e^-\mid C_e^+,C_e^-)
=1-\rho_e=\TV(p_e^+,p_e^-).
\label{eq:maximal-properties}
\end{equation}
MS-CGE preserves both marginals and attains the minimum immediate switch probability.
For coordinate $j$ and replication $i$, let $L_{j,i}^{\mathrm{MS},\pm}$ be the trajectory losses obtained using MS-CGE at every matched event. Define partial and full gradient estimates of MS-CGE as
\begin{equation*}
\widehat G_{j,n}^{\mathrm{MS}}(\theta,c)
:=\frac{1}{n}\sum_{i=1}^n
\frac{L_{j,i}^{\mathrm{MS},+}-L_{j,i}^{\mathrm{MS},-}}{2c},\qquad
\widehat{\boldsymbol g}_{n}^{\mathrm{MS}}(\theta,c)
:=\sum_{j=1}^d\widehat G_{j,n}^{\mathrm{MS}}(\theta,c)e_j.
\end{equation*}
Here, replications are independent within each coordinate, and coupling occurs only between the $+$ and $-$ trajectories in each pair. A full MS-CGE gradient uses $dn$ pairs of trajectories, with no additional score evaluations for behavioral coupling.

\subsection{CGE with First-Switch Conditioning}
\label{subsec:fs-cge}

Maximal coupling reduces how often the paired actions disagree, but each remaining switch may still produce a large finite-difference contribution. CGE with first-switch conditioning (FS-CGE) addresses this variation by averaging over the no-switch and first-switch cases. It constructs a common-action trajectory pair, called the \emph{backbone}, and simulates at most one continuation in which a switch occurs. The impact-aware extension further selects the action switch based on its predicted loss impact. We provide the details of FS-CGE in Algorithm \ref{alg:fs-cge}.

\paragraph{First-switch decomposition.}
At each backbone event $e$, evaluate $p_e^\pm$, $h_e$, and $\rho_e$ at the current contexts as in Definition~\ref{def:external-maximal}. If $\rho_e>0$, draw a common action from $h_e/\rho_e$ and apply it to both systems. Their continuous states may differ even though their actions agree. With $\delta_e=1-\rho_e$, define
\begin{equation}
W_0=1,\qquad w_e=W_{e-1}\delta_e,\qquad
W_e=W_{e-1}\rho_e,\qquad
S=\sum_{e=1}^m w_e=1-W_m.
\label{eq:first-switch-weights}
\end{equation}
The weights $w_e$ and $W_m$ represent first-switch and no-switch cases. In the first case, the weight for event $e$ includes agreement at every earlier event. Averaging these backbone-dependent weights gives the corresponding MS-CGE probabilities. Let $\Delta_0$ be the full trajectory loss difference on a completed backbone. If an overlap $\rho_e$ is zero, stop before drawing a common action, set subsequent weights to zero, and take $W_m=0$ and $\Delta_0=0$. 

For the first-switch case, a residual continuation starts from the two saved pre-action states at an event with $w_e>0$. It draws $A_e^+$ and $A_e^-$ independently from
\[
\frac{p_e^+-h_e}{\delta_e}
\quad\text{and}\quad
\frac{p_e^--h_e}{\delta_e},
\]
respectively, then uses the MS-CGE pair oracle through the remaining events. The two residual distributions have disjoint supports, so this draw forces the first switch at $e$. Denote the resulting full trajectory loss difference by $\Delta_e$, including losses accumulated before $e$. The continuation reuses the physical primitives, with fresh behavioral draws independent of backbone actions. Specifically, for $S>0$, FS-CGE selects $R$ conditionally on the backbone with probabilities $q_e^0=w_e/S$ and simulates its residual continuation.
Algorithm~\ref{alg:fs-cge} implements the selection with weighted reservoir sampling and one checkpoint. At event $e$, it replaces the checkpoint with probability $w_e/S_e$, where $S_e=\sum_{r\le e}w_r$. This checkpoint stores both restartable states, including histories and accumulated losses, $p_R^\pm$, and physical primitives. 

Combining the branches with their no-switch and first-switch weights forms FS-CGE. One coordinate contribution is then constructed as
\[
G_j^{\mathrm{FS}}(\theta,c)
=\frac{W_m\Delta_0+S\Delta_R}{2c},
\]
and the corresponding gradient estimates are
\begin{equation}
\begin{aligned}
\widehat G_{j,n}^{\mathrm{FS}}(\theta,c)
=\frac1n\sum_{i=1}^nG_{j,i}^{\mathrm{FS}}(\theta,c),\qquad
\widehat{\boldsymbol g}_n^{\mathrm{FS}}(\theta,c)
=\sum_{j=1}^d\widehat G_{j,n}^{\mathrm{FS}}(\theta,c)e_j.
\end{aligned}
\label{eq:fs-estimator}
\end{equation}
When $S=0$, set $G_j^{\mathrm{FS}}(\theta,c)=\Delta_0/(2c)$ and omit the residual continuation.
This averaging removes the Bernoulli variation in whether a switch occurs, while its location remains random.
We will see later that this weighted branch mixture enjoys marginal preservation and uniformly bounded variance (see Section~\ref{subsec:fs-theory} below).
Notably, the impacts on the full trajectory loss may vary significantly across events, which motivates the first-switch event selection that also accounts for this effect.

\begin{algorithm}[!t]
\caption{One FS-CGE contribution for coordinate $j$}
\label{alg:fs-cge}
\renewcommand{\baselinestretch}{1}\fontsize{10}{12}\selectfont
\begin{algorithmic}[1]
\Require Policy $\theta$, coordinate $j$, perturbation radius $c>0$, and simulation horizon $T$.
\State Draw shared physical randomness $\boldsymbol\xi$ and initialize both simulators at $\theta\pm ce_j$.
\State $W\gets1$, $S\gets0$, $R\gets0$, $\Delta_0\gets0$, $\Delta_R\gets0$.
\For{each behavioral event $e=1,\ldots,m$}
    \Comment{Build the common-action backbone.}
    \State Advance both simulations to event $e$ using $\boldsymbol\xi$ and evaluate $p_e^\pm$, $h_e$, and $\rho_e$.
    \State $w_e\gets W(1-\rho_e)$, $S\gets S+w_e$.
    \If{$w_e>0$}
        \Comment{Select the checkpoint.}
        \State With probability $w_e/S$, set $R\gets e$ and overwrite the saved data with both full pre-action states and $p_e^\pm$.
    \EndIf
    \State $W\gets W\rho_e$. Stop the loop if $\rho_e=0$.
    \State Draw a common action $A_e\sim h_e/\rho_e$ and apply it to both simulations.
\EndFor
\State Set $W_m\gets W$. If $W_m>0$, continue to $T$ and set $\Delta_0$ to the total loss difference.
\If{$S>0$}
    \Comment{Restore the switch at $R$ and continue.}
    \State Restore the saved states and probability vectors $p_R^\pm$.
    \State Compute $h_R,\rho_R$ from $p_R^\pm$ and independently draw $A_R^\pm\sim(p_R^\pm-h_R)/(1-\rho_R)$.
    \State Apply $A_R^\pm$ and continue to $T$ using MS-CGE and the same $\boldsymbol\xi$, recomputing contexts and probabilities.
    \State Set $\Delta_R$ to the total loss difference, including the saved losses.
\EndIf
\State \Return $G_j^{\mathrm{FS}}(\theta,c)=(W_m\Delta_0+S\Delta_R)/(2c)$.
\end{algorithmic}
\end{algorithm}

\paragraph{Impact-aware first-switch selection.}
FS-CGE selects the first-switch event using the first-switch weights with probability $q_e^0=w_e/S$. Here, we introduce CGE with impact-aware first-switch conditioning (IA-FS-CGE) that further accounts for the predicted magnitude of the full trajectory loss difference $\Delta_e$. It uses the same backbone and at most one residual continuation but with adjusted selection probabilities and branch weights.

We assume access to a predictor of the squared full loss difference $\Delta_e^2$, fitted independently of the evaluation replications and held fixed during gradient estimation. Its inputs include the policy, perturbation radius, event, profile, and state information available before the common action at event $e$. The predicted second moment accounts for both the variance and squared mean of the loss difference. We take its square root to obtain a predicted root mean square (RMS), which is then normalized by an shared scale across events at the current policy. Clipping the resulting multiplier to $[\underline u,\overline u]$ gives $u_e$, where $0<\underline u\le1\le\overline u<\infty$.

Let $\mathcal B$ contain the physical primitives and sampled backbone. The product $w_eu_e$ is the sampling weight for event $e$. Its cumulative sum is $U_e=\sum_{r\le e}w_ru_r$, so $U_m$ is the total weight used to normalize the selection probabilities. For $S>0$, define
\begin{equation}
\widehat q_e=\frac{w_eu_e}{U_m},\qquad
q_e^\gamma=(1-\gamma)\frac{w_e}{S}+\gamma\widehat q_e,
\qquad 0\le\gamma\le1.
\label{eq:ia-fs-proposal}
\end{equation}
Here $\widehat q_e$ combines first-switch weight with predicted loss impact, while $\gamma$ controls how strongly the predicted impact changes the final probability $q_e^\gamma$. Setting $\gamma=0$ recovers FS-CGE, and $\gamma=1$ uses the predicted weights alone. The clipping of $u_e$ keeps the selection probability of each event with a positive weight at least a fixed positive fraction of its FS-CGE value.

For $S>0$, let $q(\cdot\mid\mathcal B)$ be any probability vector on the positive-weight events, with $q_e>0$ for each such event. We set $q_e=q_e^\gamma$ using the frozen predictor and draw the first-switch event $R$ from this distribution. We then simulate the residual continuation from event $R$ as in FS-CGE and observe its actual full loss difference $\Delta_R$. The single-replication contribution of IA-FS-CGE is
\begin{equation}
G_j^{\mathrm{IA}}(\theta,c;q)
=\frac{W_m\Delta_0+(w_R/q_R)\Delta_R}{2c}.
\label{eq:ia-fs-contribution}
\end{equation}
The term $W_m\Delta_0$ is the weighted no-switch contribution. In the residual term, the factor $w_R/q_R$ compensates for sampling event $R$ with probability $q_R$ to restore its first-switch weight $w_R$ in expectation. The probabilities are determined before the residual continuation is generated, so this correction preserves the central-difference mean even when the predictions are inaccurate. When $S=0$, set $G_j^{\mathrm{IA}}(\theta,c;q)=\Delta_0/(2c)$. Averaging evaluation replications that are independent conditional on the frozen predictor gives
\begin{equation}
\widehat G_{j,n}^{\mathrm{IA}}(\theta,c;q)
=\frac1n\sum_{i=1}^nG_{j,i}^{\mathrm{IA}}(\theta,c;q),\qquad
\widehat{\boldsymbol g}_n^{\mathrm{IA}}(\theta,c;q)
=\sum_{j=1}^d\widehat G_{j,n}^{\mathrm{IA}}(\theta,c;q)e_j.
\label{eq:ia-fs-estimator}
\end{equation}

The probabilities $q_e^\gamma$ can be implemented with the same single checkpoint used by FS-CGE. Specifically, draw $H\sim\mathrm{Bernoulli}(\gamma)$ once per replication. If $H=0$, use reservoir weights $w_e$, and if $H=1$, use $w_eu_e$. Equivalently, write $v_e=w_e[(1-H)+Hu_e]$ and $V_e=\sum_{r\le e}v_r$. At event $e$, replace the checkpoint with probability $v_e/V_e$ when $v_e>0$, and save $u_e$ with the restart data. Accumulate $S$ and $U_m$ along the backbone, and averaging the two possible weight rules over $H$ gives the selection probability $q_e^\gamma$. Therefore, the final correction uses this marginal mixture probability:
\begin{equation}
\frac{w_R}{q_R^\gamma}
=\left(\frac{1-\gamma}{S}+\frac{\gamma u_R}{U_m}\right)^{-1}.
\label{eq:ia-fs-mixture-weight}
\end{equation}
This factor equals $S$ when $\gamma=0$ and $U_m/u_R$ when $\gamma=1$. For intermediate values of $\gamma$, both terms in the denominator are used regardless of the realized value of $H$. Thus, $H$ selects the reservoir weights, while $q_R^\gamma$ determines the correction.

For both FS-CGE and IA-FS-CGE, building the backbone requires at most $2m$ score evaluations, and reusing $p_R^\pm$ leaves at most $2(m-R)$ score evaluations for the residual continuation. The total is then at most $4m-2$ when a residual is needed. Besides, predictor training and evaluation for IA-FS-CGE may incur additional work beyond the online score evaluations in practice.

\section{Analysis of CGEs with Maximal Coupling}\label{sec:analysis_PGE}
The analysis concerns the objective $J$ under the score--softmax construction \eqref{eq:score-softmax-kernel}. We first relate action switches to the variance of MS-CGE, and then show how first-switch conditioning and residual allocation affect that variance and the cost of optimization.

\begin{definition}
For $\theta\in\Theta$, coordinate $j$, and $c>0$, an \emph{action switch} occurs at event $e$ when $A_{e,j}^{+}\ne A_{e,j}^{-}$. Define the action switch count in a maximally coupled trajectory pair as
\[
N_{\mathrm{as},j}^{\mathrm{MS}}(\theta,c)
:=\sum_{e=1}^m\ind\{A_{e,j}^{+}\ne A_{e,j}^{-}\}.
\]
\end{definition}

The next assumption bounds the loss effects of policy perturbations and action switches.

\begin{assumption}\label{assump:switch_dominance}
There exist nonnegative constants $L_f<\infty$ and $C_{\mathrm{as}}<\infty$ such that, uniformly over $\theta\in\Theta$, coordinates $j$, and $0<c\le c_0$, the loss difference of a two-sided trajectory pair with shared physical randomness satisfies
\[
\left|L_{j}^{r,+}-L_{j}^{r,-}\right|
\le L_f c+C_{\mathrm{as}}N_{\mathrm{as},j}^{r}(\theta,c),
\]
where $N_{\mathrm{as},j}^{r}$ counts its action switches.
\end{assumption}
The bound applies to histories supported by I-CGE and MS-CGE and to the backbone and residual continuations of FS-CGE and IA-FS-CGE. It holds when loss is uniformly Lipschitz in the policy with actions fixed, with a common deterministic constant over supported primitives and histories, and each action change has a uniformly bounded downstream effect. We also require the following assumption.

\begin{assumption}\label{assump:switch_sparsity}
There exist constants $C_{\mathrm{sp}}<\infty$, $\kappa>0$, and $c_0>0$ such that, uniformly over $\theta\in\Theta$, coordinates $j$, and $0<c\le c_0$,
\[
\E\!\left[N_{\mathrm{as},j}^{\mathrm{MS}}(\theta,c)\right] \le C_{\mathrm{sp}}c^{\kappa}.
\]
\end{assumption}
The constant $C_{\mathrm{sp}}$ sets the scale of the bound, and $\kappa$ describes how quickly the expected switch count decreases as the perturbation radius shrinks.
Assumptions~\ref{assump:smoothness} and~\ref{assump:switch_dominance} apply throughout, and we will provide sufficient switch-sparsity conditions for Assumption \ref{assump:switch_sparsity} below.

\subsection{Bias and Variance of MS-CGE}

Exact marginality preserves each standalone trajectory law, even after the two histories diverge. Thus $\E[\widehat G_{j,n}^{\mathrm{MS}}(\theta,c)]=D_{j,c}J(\theta)$, and Lemma~\ref{lemma:pge-mean-preservation} in Appendix~\ref{app:proof-pge-mean} gives the event-wise induction. The next theorem combines finite-difference bias with switch-induced variance.

\begin{theorem}
\label{thm:pge-main}
Under Assumptions~\ref{assump:smoothness}, \ref{assump:switch_dominance}, and \ref{assump:switch_sparsity}, each coordinate of MS-CGE satisfies
\begin{equation}
\begin{aligned}
\left|\E\!\left[\widehat G_{j,n}^{\mathrm{MS}}(\theta,c)\right]-\partial_jJ(\theta)\right|
\le \frac{L_Jc}{2},\qquad
\Var\!\left(\widehat G_{j,n}^{\mathrm{MS}}(\theta,c)\right)
 \le \frac{L_f^2}{2n}
+\frac{C_{\mathrm{as}}^2mC_{\mathrm{sp}}}{2n}c^{\kappa-2}.
\label{eq:pge-component-variance}
\end{aligned}
\end{equation}
For the full MS-CGE vector,
\begin{equation}
\begin{aligned}
\left\|\E[\widehat{\boldsymbol g}_{n}^{\mathrm{MS}}(\theta,c)]-\nabla J(\theta)\right\|^2
&\le \frac{dL_J^2c^2}{4},\\
\E\!\left[\left\|\widehat{\boldsymbol g}_{n}^{\mathrm{MS}}(\theta,c)-\nabla J(\theta)\right\|^2\right]
&\le \frac{dL_J^2c^2}{4}
+\frac{dL_f^2}{2n}
+\frac{dC_{\mathrm{as}}^2mC_{\mathrm{sp}}}{2n}c^{\kappa-2}.
\end{aligned}
\label{eq:pge-vector-mse}
\end{equation}
\end{theorem}

In \eqref{eq:pge-vector-mse}, the terms in the MSE bound come from finite-difference bias, smooth path differences, and action switches. The radius controls the tradeoff between the first and third terms.

For $\kappa=1$, the radius is of order $n^{-1/3}$ and the MSE bound is $O(dn^{-2/3})$. Once $\kappa\ge2$, the switch term is bounded as $c$ decreases, and a radius of order $n^{-1/2}$ yields $O(d/n)$ MSE. We provide in Appendix~\ref{app:proof-pge-main} the proof of Theorem~\ref{thm:pge-main} and the general radius schedule in Corollary~\ref{cor:pge-mse-rate}. Under the Lipschitz-gradient assumption, these schedules balance the stated bias and variance bounds. 

\subsection{Switch Sparsity Before the First Switch}
\label{subsec:switch-stability}

To verify the switch-sparsity condition from the model, it suffices to control trajectories up to their first switch. Over a fixed horizon, a small first-switch probability bounds the expected total action switches even when the trajectories subsequently diverge.

Let $\sigma_j=\min\{e:A_{e,j}^+\ne A_{e,j}^-\}$ be the first-switch event, with $\sigma_j=\infty$ if no switch occurs. Let $\mathcal F_{e-1}$ contain the physical primitives and paired histories before event $e$. The event $\{\sigma_j\ge e\}$ is $\mathcal F_{e-1}$-measurable, and $\{\sigma_j\ge e\}$ means that all actions before $e$ agree. The stability condition concerns these common-action prefixes and permits divergence after a switch.

\begin{assumption}
\label{assump:preswitch}
For some $0<\alpha\le1$ and finite nonnegative deterministic constants $M_1,\ldots,M_m$, uniformly over $\theta$, $j$, and sufficiently small $0<c\leq c_0$,
\begin{equation*}
\delta_{e,j}:=\TV(p_{e,j}^+,p_{e,j}^-)
\le M_e c^\alpha
\quad\text{on }\{\sigma_j\ge e\},\qquad e=1,\ldots,m.
\end{equation*}
\end{assumption}
Here $M_e$ bounds event-specific sensitivity and $\alpha$ describes local regularity of the behavioral probabilities. This pre-switch bound implies switch sparsity, as shown below.
\begin{proposition}
\label{prop:switchsparse}
Under Assumption~\ref{assump:preswitch},
\[
\E[N_{\mathrm{as},j}^{\mathrm{MS}}(\theta,c)]
\le c^\alpha\sum_{e=1}^m(m-e+1)M_e.
\]
Thus, Assumption~\ref{assump:switch_sparsity} holds with $\kappa=\alpha$.
\end{proposition}

Exact maximal coupling gives $\Prob(\sigma_j=e)=\E[\ind\{\sigma_j\ge e\}\delta_{e,j}]\le M_ec^\alpha$. After a first switch at $e$, there are at most $m-e+1$ switches in total. This gives Proposition~\ref{prop:switchsparse}, whose proof is provided in Appendix~\ref{app:proof-switchsparse}. In fact, Assumption~\ref{assump:preswitch} bounds the disagreement probability at each event on every supported common-action prefix, which imposes stronger requirements than Assumption~\ref{assump:switch_sparsity}.

A sufficient condition can be checked from the simulator. Let $\widetilde X_e^\pm$ include the pre-action state, context history, and event activity, and suppose we can define the distance $D_e=\|\widetilde X_e^+-\widetilde X_e^-\|$. On common-action prefixes, suppose finite nonnegative constants satisfy
\begin{equation}
D_1\le\beta_0c,\qquad D_{e+1}\le r_eD_e+\beta_ec,\qquad
\|s_e(C_e^+)-s_e(C_e^-)\|_\infty\le L_sD_e^\alpha.
\label{eq:preswitch-primitives}
\end{equation}
Define $K_1=\beta_0$ and $K_{e+1}=r_eK_e+\beta_e$. Induction using \eqref{eq:preswitch-primitives} then gives $D_e\le K_ec$. The softmax bound
\[
\TV(\operatorname{softmax}_\tau(x),\operatorname{softmax}_\tau(y))\le\|x-y\|_\infty/\tau
\]
then yields $M_e=L_sK_e^\alpha/\tau$. The conditions concern pre-switch dynamics and permit post-switch divergence.

For $\alpha=1$, Proposition~\ref{prop:switchsparse} gives an $O(c)$ switch bound. The next result shows that this order is sharp under nonzero first-order behavioral sensitivity.

\begin{proposition}
\label{prop:first-order-switch}
Fix a policy \(\theta\) and coordinate \(j\). With the physical primitives \(\boldsymbol\xi\) held fixed, suppose that, almost surely in \(\boldsymbol\xi\), the partial derivative \(\partial_j P_{b,1}(a\mid C_1^\theta)\) exists for every \(a\in\mathcal A\), and \[
\mathbb P \bigg(\sum_{a\in\mathcal A}
\left|\partial_j P_{b,1}(a\mid C_1^\theta)\right|>0\bigg)>0.
\]
Then, \(\mathbb E[N_{\mathrm{as},j}^{\mathrm{MS}}(\theta,c)]=\Omega(c)\). Consequently, an upper bound of the form \(O(c^\kappa)\) with \(\kappa>1\) cannot hold at this policy and coordinate.
\end{proposition}

The proof is given in Appendix~\ref{app:proof-first-order}. Under the sensitivity condition in Proposition~\ref{prop:first-order-switch} and locally Lipschitz pre-switch stability, Proposition~\ref{prop:switchsparse} gives a matching $O(c)$ upper bound, so the expected switch count is $\Theta(c)$. Higher exponents require first-order insensitivity, which smoothness alone does not imply.
Appendix~\ref{app:proof-first-order} also provides a binary example in which maximal coupling minimizes the variance of the paired finite difference, while the variance grows as \(c\) decreases, illustrating the motivation for first-switch conditioning.

\subsection{Bias and Variance of FS-CGE}
\label{subsec:fs-theory}

For one coordinate, let $G_j^{\mathrm{MS}}(\theta,c)=(L_j^{\mathrm{MS},+}-L_j^{\mathrm{MS},-})/(2c)$. The next result gives the variance reduction per replication from first-switch conditioning.

\begin{theorem}
\label{thm:fs-conditioning}
Under Assumption~\ref{assump:switch_dominance} and the construction in Section~\ref{subsec:fs-cge}, FS-CGE has the same central-difference mean as MS-CGE, $\E[G_j^{\mathrm{FS}}(\theta,c)]=D_{j,c}J(\theta)$, and
\begin{equation}
\Var(G_j^{\mathrm{MS}}(\theta,c))-\Var(G_j^{\mathrm{FS}}(\theta,c))
=\frac{\E[S(1-S)(\Delta_R-\Delta_0)^2]}{4c^2}\ge0.
\label{eq:fs-variance-identity}
\end{equation}
\end{theorem}

FS-CGE removes the Bernoulli variance $S(1-S)$ of choosing between the two cases, weighted by their squared loss difference. The estimator still samples the physical inputs, backbone, switch location and continuation. Its advantage is therefore small when the two cases have similar losses or one already has probability close to one.

Let $M$ collect the event-wise sensitivity bounds and let $B$ bound the magnitude of a residual loss difference. Define $M=\sum_{e=1}^m M_e$, $B=L_fc_0+C_{\mathrm{as}}m$, and $C_{\mathrm{FS}}=\frac{L_f+MB}{2}$.
Then, the following theorem provides the bias and variance bounds of FS-CGE.

\begin{theorem}
\label{thm:fs-rate}
Under Assumptions~\ref{assump:smoothness}, \ref{assump:switch_dominance} and \ref{assump:preswitch} with $\alpha=1$, each coordinate satisfies
\[
|\E[\widehat G_{j,n}^{\mathrm{FS}}(\theta,c)]-\partial_jJ(\theta)|\le\frac{L_Jc}{2},
\qquad
\Var(\widehat G_{j,n}^{\mathrm{FS}}(\theta,c))\le\frac{C_{\mathrm{FS}}^2}{n}.
\]
Consequently,
\begin{equation}
\E\|\widehat{\boldsymbol g}_n^{\mathrm{FS}}(\theta,c)-\nabla J(\theta)\|^2
\le\frac{dL_J^2c^2}{4}+\frac{dC_{\mathrm{FS}}^2}{n}.
\label{eq:fs-mse}
\end{equation}
\end{theorem}
The variance bound follows from $S\le Mc$, $|\Delta_0|\le L_fc$, and $|\Delta_R|\le B$. Thus, both the no-switch and first-switch terms are $O(c)$ in a FS-CGE contribution. Moreover, choosing $c_n=\Theta(n^{-1/2})$ leads to an $O(d/n)$ MSE. For \(\alpha<1\), the pre-switch condition does not guarantee a uniform variance bound.

The residual allocation in IA-FS-CGE affects variance through the loss differences of the sampled branches. Given the physical inputs and sampled backbone in $\mathcal B$, let $\mathcal M_e=\E[\Delta_e^2\mid\mathcal B]$.
This second moment includes the within-branch variance and squared conditional mean, so $\sqrt{\mathcal M_e}$ is the conditional RMS of the full loss difference for a first switch at event $e$.
For $S>0$, the FS-CGE probabilities are $q_e^0=w_e/S$. The oracle benchmark uses $q_e^\star=w_e\sqrt{\mathcal M_e}/Z$ when $Z=\sum_e w_e\sqrt{\mathcal M_e}>0$. 

\begin{proposition}
\label{prop:ia-fs-allocation}
Under Assumption~\ref{assump:switch_dominance} and the construction in Section~\ref{subsec:fs-cge}, fix the predictor used to obtain the selection probabilities. Every admissible $q$ preserves $\E[G_j^{\mathrm{IA}}(\theta,c;q)]=D_{j,c}J(\theta)$. The oracle variance gain is
\begin{equation}
\Var(G_j^{\mathrm{FS}}(\theta,c))-\Var(G_j^{\mathrm{IA}}(\theta,c;q^\star))
=\frac{\E[S^2\Var_{e\sim q^0}(\sqrt{\mathcal M_e})]}{4c^2}\ge0.
\label{eq:ia-fs-oracle-gain}
\end{equation}
For any admissible $q$ with $\Var(G_j^{\mathrm{IA}}(\theta,c;q))<\infty$, define $\chi^2(q^\star\Vert q)=\sum_{e:q_e>0}\frac{(q_e^\star-q_e)^2}{q_e}$. Then,
\begin{equation}
\begin{aligned}
\Var(G_j^{\mathrm{IA}}(\theta,c;q))-\Var(G_j^{\mathrm{IA}}(\theta,c;q^\star))
=\frac{\E[Z^2\chi^2(q^\star\Vert q)]}{4c^2}.
\end{aligned}
\label{eq:ia-fs-allocation-penalty}
\end{equation}
\end{proposition}

In \eqref{eq:ia-fs-oracle-gain}, $\Var_{e\sim q^0}$ measures variation across switch locations on a fixed backbone, so the oracle gain comes from differences in their RMS loss impacts. The predictor uses information available before each action, while the oracle uses conditional moments given the full backbone. Prediction errors affect variance through the resulting selection probabilities. In \eqref{eq:ia-fs-allocation-penalty}, the $\chi^2$ term measures their mismatch with $q^\star$. Underestimating a branch with a large RMS loss impact can assign it too little probability and increase this penalty. IA-FS-CGE reduces variance relative to FS-CGE when the penalty is smaller than the oracle gain. When RMS impacts are similar across branches, the oracle gain itself is small.

Clipping $u_e$ to $[\underline u,\overline u]$ limits the importance weights. Write $\Lambda=\overline u/\underline u$ for the ratio of the clipping bounds. Since $u_e\ge\underline u$ and $U_m\le\overline u S$, we have $q_e^\gamma\ge q_e^0/\Lambda$ and hence $w_R/q_R^\gamma\le\Lambda S$ for every $\gamma\in[0,1]$.

\begin{corollary}
\label{cor:ia-fs-rate}
Under the assumptions of Theorem~\ref{thm:fs-rate}, define $C_{\mathrm{IA}}=(L_f+\Lambda MB)/2$.
Then,
$|G_j^{\mathrm{IA}}(\theta,c;q^\gamma)|\le C_{\mathrm{IA}}$ for every $\gamma\in[0,1]$, and
\begin{equation}
\E\|\widehat{\boldsymbol g}_n^{\mathrm{IA}}(\theta,c;q^\gamma)-\nabla J(\theta)\|^2
\le\frac{dL_J^2c^2}{4}+\frac{dC_{\mathrm{IA}}^2}{n},
\label{eq:ia-fs-mse}
\end{equation}
where the expectation includes any independent randomness used to obtain the predictor. With $c_n=\Theta(n^{-1/2})$, the MSE is $O(d/n)$.
\end{corollary}

The factor $\Lambda$ changes the variance bound's constant while preserving the FS-CGE MSE rate. Online efficiency also depends on continuation length, because earlier first switches require longer residual simulations. As a result, the experiments compare the variance--cost product using expected score-evaluation cost.

\subsection{Projected Stochastic Approximation}

We use the estimators in projected stochastic approximation. At iteration $k$, define
\[
\widehat{\boldsymbol g}_k^r
=\sum_{j=1}^d\left\{\frac1n\sum_{i=1}^nG_{k,j,i}^r\right\}e_j
\]
for $r\in\{\mathrm{I},\mathrm{MS},\mathrm{FS},\mathrm{IA}\}$, and $\theta_{k,j}^{\pm}=\theta_k\pm ce_j$.
Here I-CGE and MS-CGE use the paired loss difference divided by $2c$, while FS-CGE and its IA extension use their weighted contributions from Section~\ref{subsec:fs-cge}. Each method targets the same central difference with shared physical primitives. The update of projected stochastic approximation is
\begin{equation*}
\theta_{k+1}=\Pi_{\Theta}\!\left(\theta_k-\eta\widehat{\boldsymbol g}_k^r\right).
\end{equation*}
Let
\[
\mathcal G_\eta(\theta)
:=\eta^{-1}\left[\theta-\Pi_\Theta\!\left(\theta-\eta\nabla J(\theta)\right)\right]
\]
denote the gradient mapping, which is commonly used to measure stationarity of solutions to constrained nonconvex optimization.

For the optimization bounds, the relevant assumptions hold uniformly along the iterates, $J\ge J_{\inf}$ on $\Theta$, and $0<\eta\le1/(2L_J)$. Both FS-CGE and IA-FS-CGE use the pathwise Lipschitz condition in Assumption~\ref{assump:preswitch} with $\alpha=1$. Under these conditions, the next theorem bounds the replications needed to attain a prescribed stationarity tolerance under fixed step sizes.

\begin{theorem}
\label{thm:sample_complexity}
Assumptions~\ref{assump:smoothness} and \ref{assump:switch_dominance} apply throughout. The MS-CGE bound also uses Assumption~\ref{assump:switch_sparsity}, and the FS-CGE bound uses Assumption~\ref{assump:preswitch} with $\alpha=1$. IA-FS-CGE has the same replication order as FS-CGE, with the variance constant from Corollary~\ref{cor:ia-fs-rate}.
For sufficiently small $\epsilon>0$, the total numbers of coordinate replications sufficient for $\min_{0\le k<K}\E\|\mathcal G_\eta(\theta_k)\|^2\le\epsilon$ are
\[
\begin{array}{ll}
\text{I-CGE}:O(d^3\epsilon^{-3}),\qquad
\text{MS-CGE}:
\begin{cases}
O\!\left(d^{3-\kappa/2}\epsilon^{-3+\kappa/2}\right),&0<\kappa<2,\\
O(d^2\epsilon^{-2}),&\kappa\ge2,
\end{cases}\\
\text{FS-CGE and IA-FS-CGE}:O(d^2\epsilon^{-2}).
\end{array}
\]
\end{theorem}
Each I-CGE or MS-CGE replication is a two-sided trajectory pair with at most $2m$ score evaluations. An FS-CGE or IA-FS-CGE replication includes a backbone and at most one residual continuation, requiring at most $4m$ evaluations. Therefore, counting score evaluations adds a factor of $m$ to these cost bounds, and the variance constants can also depend on the horizon. For $\kappa=1$, MS-CGE requires $O(d^{5/2}\epsilon^{-5/2})$ coordinate replications. The stronger pathwise condition and uniform variance bound give FS-CGE and IA-FS-CGE the $O(d^2\epsilon^{-2})$ bound.

The paired construction also applies to finite differences along other directions. For $u\in\R^d$ with $c\|u\|\le c_0$, replace $e_j$ by $u$ to estimate
\[
D_{u,c}J(\theta)=\frac{J(\theta+cu)-J(\theta-cu)}{2c}.
\]
Conditional on $u$, maximal coupling and first-switch conditioning preserve this directional mean and can reduce simulation variance under corresponding directional stability conditions. These scalar estimates can be used in simultaneous perturbation stochastic approximation (SPSA) \citep{spall1992multivariate} and other randomized gradient estimators \citep{berahas2022theoretical}. Random directions introduce variation beyond simulation noise. The direction distribution and reconstruction weights determine the resulting gradient estimator's bias and variance, so its stochastic approximation bounds require a separate directional analysis.

\section{Numerical Simulations}\label{sec:simulation}

We compare the four estimators in a synthesized service-feedback model and an LLM-embedded software-onboarding model. Within each setting, all methods share the behavioral law, loss, and physical randomness. We use I, MS, FS, and IA to represent the four estimators. For optimization, we report the mean of $J(\theta)-J(\theta^\star)$ over runs, where $\theta^\star$ minimizes $J$ for each model. Online budgets exclude training costs and assume predictor reuse. Detailed problem formulation and parameter selection are given in Appendix~\ref{app:numerical-settings}.

\subsection{Synthesized Experiments Without an LLM}
\label{subsec:controlled-maximal-test}

We study a service-feedback problem over four interactions ($m=4$) with one binary satisfaction response per interaction. The effort policy $\theta\in[0.1,1.9]$ balances a quadratic effort cost against dissatisfaction. A persistent customer type and response feedback link successive interactions, allowing an action switch to affect later losses.

\begin{figure}[!htbp]
    \centering
    \includegraphics[width=0.90\textwidth]{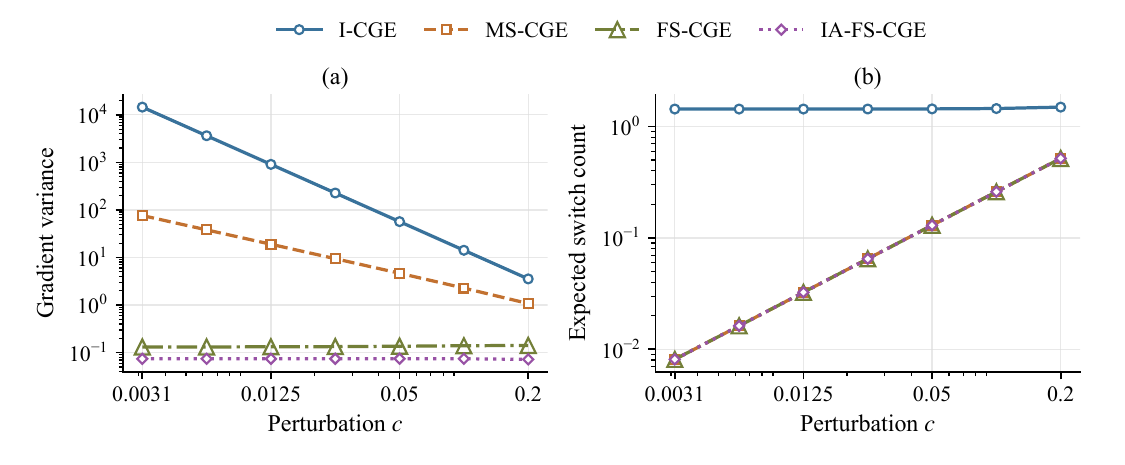}
    \caption{Gradient variance (a) and switch diagnostics (b) in the synthesized model.}
    \label{fig:controlled-mechanism}
\end{figure}

Figure~\ref{fig:controlled-mechanism} compares variance and switch counts at $\theta=0.2$. As $c$ decreases, the variance of I rises sharply while its switch frequency remains almost unchanged. With MS, switches become less frequent as $c$ decreases, although the surviving switched-path loss differences still lead to increasing variance. These trends agree with the $c^{-2}$ and $c^{-1}$ variance bounds. FS variance remains nearly constant across the tested radii, and IA maintains a lower variance throughout. Here, $N_{\mathrm{as}}^{\mathrm{branch}}$ counts switches in the sampled residual continuation, including the forced first switch, and is zero when no continuation is generated. The weighted counts $S N_{\mathrm{as}}^{\mathrm{branch}}$ and $(w_R/q_R)N_{\mathrm{as}}^{\mathrm{branch}}$ of FS and IA retain the same expectation as MS. Both FS and IA reduce the variance of switched-path contributions.

Enumeration verifies the common central-difference mean and MS--FS variance identity in Theorem~\ref{thm:fs-conditioning}, so all four estimators have the same finite-difference bias. The IA improvement persists when continuation work is counted. To assess efficiency over the policy range, we compare the variance--evaluation-cost products on an equally spaced grid over \(\theta\in[0.1,1.9]\). At \(c=0.05\), IA achieves an average reduction of \(43.8\%\) relative to FS.

Figure~\ref{fig:controlled-mse} shows how these variance differences affect gradient MSE under the radius schedules in Appendix~\ref{app:controlled-tuning}. The I and MS curves decay more slowly than the nearly parallel FS and IA curves, consistent with the corresponding MSE bounds. All curves use exact moments, with IA averaging twenty independent fits per radius. IA reduces MSE by $44$--$47\%$ at the same bias compared with FS, and the nearly parallel FS and IA curves suggest similar decay rates.

\begin{figure}[!htbp]
    \centering
    \includegraphics[width=0.90\textwidth]{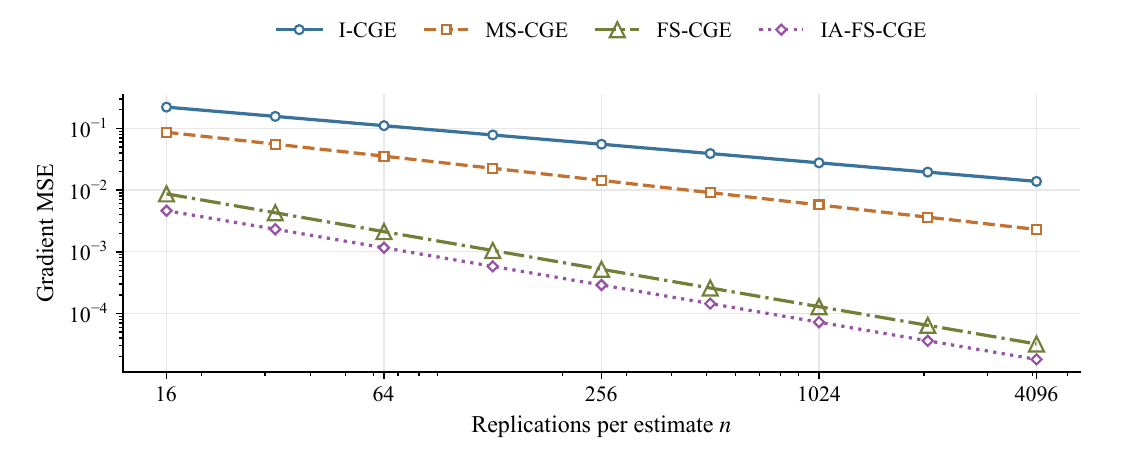}
    \caption{Exact gradient MSE in the synthesized model.}
    \label{fig:controlled-mse}
\end{figure}

\begin{figure}[!htbp]
    \centering
    \includegraphics[width=0.90\textwidth]{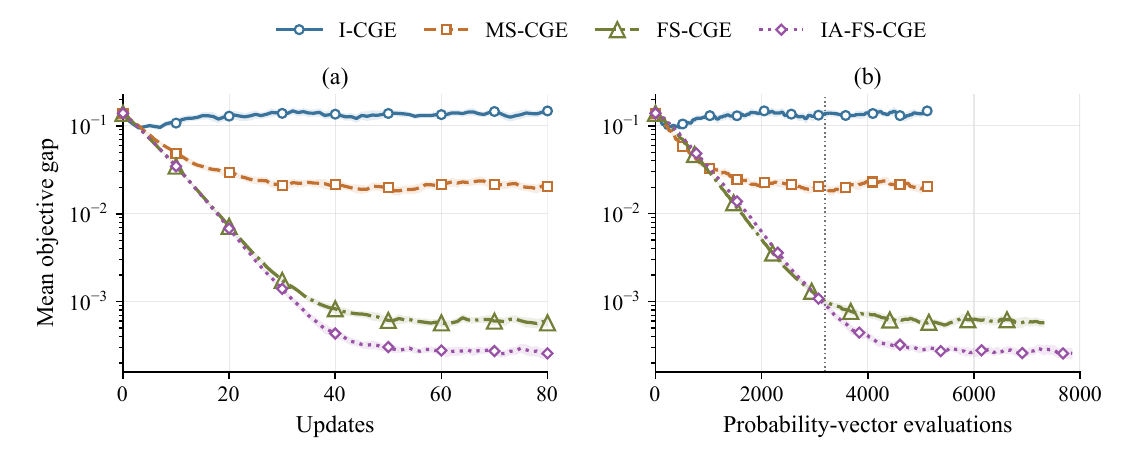}
    \caption{Synthesized optimization by updates (a) and computational cost (b).}
    \label{fig:controlled-optimization}
\end{figure}

With a common fixed step size, lower gradient noise also permits smaller optimization gaps in Figure~\ref{fig:controlled-optimization}. The runs use $\theta_0=0.2$, $c=0.05$, $n=8$ and $\eta=0.1$. Both panels show the same 80-update trajectories, and the work panel retains each run's last completed iterate. I and MS level off at larger gaps than FS and IA, consistent with the ordering of their gradient variances. IA reaches a $55.0\%$ smaller final mean gap than FS while using somewhat more evaluation calls. In panel (b), IA also maintains a lower mean gap over the later part of the common evaluation range, indicating that the improved precision can offset its additional continuation work.

\subsection{LLM-Embedded Software-Onboarding Experiments}
\label{subsec:llm-experiments}

We study a software onboarding problem over at most four invitation rounds $(m=4)$, with users drawn from two equally likely, persistent profiles. In each round, users attempt setup or defer. Successful activation ends the process, while deferral or technical failure affects later invitations. The policy $\theta\in[0.1,0.9]$ controls the probability of live assistance, balancing activation against capacity and contact costs. At each state, the attempt probability is a linear mixture of the fixed probabilities under invitations with and without assistance, with weights $\theta$ and $1-\theta$. We use Qwen3.5-2B (2B) and Qwen3-4B-Instruct-2507 (4B) to supply the action scores used by all four methods (see Appendix~\ref{app:llm-settings}).

Table~\ref{tab:llm-gradient-moments} compares exact moments at $\theta=0.8$, $c=0.025$, for one replication and an average of $n=4$. All methods share the mean and bias within each kernel. MS removes most independent-action noise, and FS further removes much of the remaining variance. IA achieves a further reduction under the 2B model, with little improvement under the 4B model.

\begin{table}[!htbp]\centering
\caption{Exact gradient moments under the LLM kernels. IA averages twenty offline fits.}\label{tab:llm-gradient-moments}
\ManuscriptTableFont
\begin{tabular}{llrrrr}\toprule
Model & Method & $\E[G]$ & Bias & $\Var(G)$ & $\Var(\widehat G_{1,4})$\\\midrule
Qwen3.5-2B & I-CGE & $0.344$ & $-2.115\!\times\!10^{-6}$ & $65.857$ & $16.464$\\
 & MS-CGE & $0.344$ & $-2.115\!\times\!10^{-6}$ & $1.615$ & $0.404$\\
 & FS-CGE & $0.344$ & $-2.115\!\times\!10^{-6}$ & $0.026416$ & $0.006604$\\
 & IA-FS-CGE & $0.344$ & $-2.115\!\times\!10^{-6}$ & $0.024439$ & $0.006110$\\
\midrule
Qwen3-4B-Instruct-2507 & I-CGE & $0.318$ & $-9.009\!\times\!10^{-9}$ & $76.659$ & $19.165$\\
 & MS-CGE & $0.318$ & $-9.009\!\times\!10^{-9}$ & $2.303$ & $0.576$\\
 & FS-CGE & $0.318$ & $-9.009\!\times\!10^{-9}$ & $0.004239$ & $0.001060$\\
 & IA-FS-CGE & $0.318$ & $-9.009\!\times\!10^{-9}$ & $0.004234$ & $0.001059$\\
\bottomrule\end{tabular}\end{table}

\begin{figure}[!htbp]
    \centering
    \includegraphics[width=0.80\textwidth]{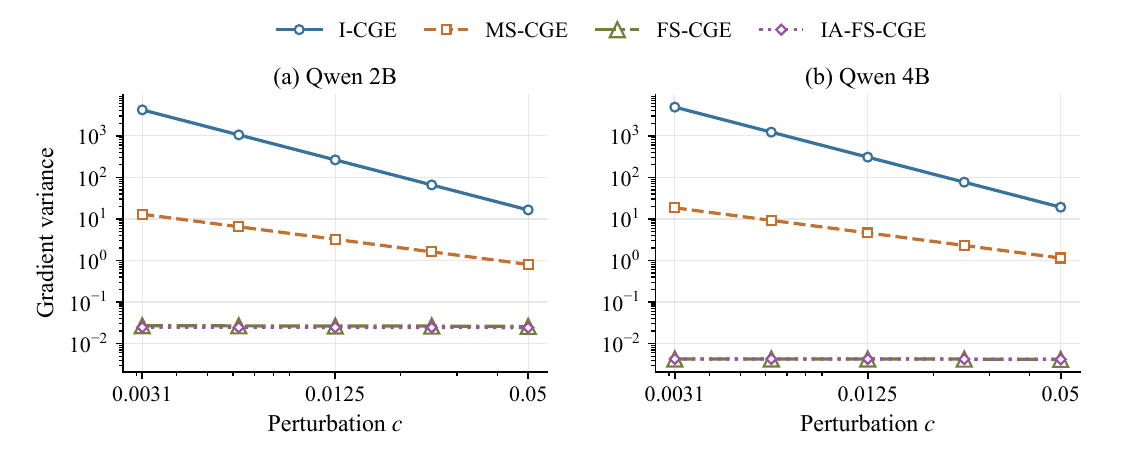}
    \caption{Gradient variance under the LLM kernels.}
    \label{fig:llm-variance-radii}
\end{figure}

The gradient variance comparisons under varying radii in Figure~\ref{fig:llm-variance-radii} reproduce the same distinction. Rare switches amplify MS noise as $c$ shrinks, while FS and IA avoid this amplification. The finite-grid trends agree with the bounds, with IA nearly coinciding with FS under the 4B model. Pre-switch stability follows from varying the assistance probability through this linear mixture while keeping the LLM scores at each state fixed.

\begin{figure}[!t]
    \centering
    \includegraphics[width=0.80\textwidth]{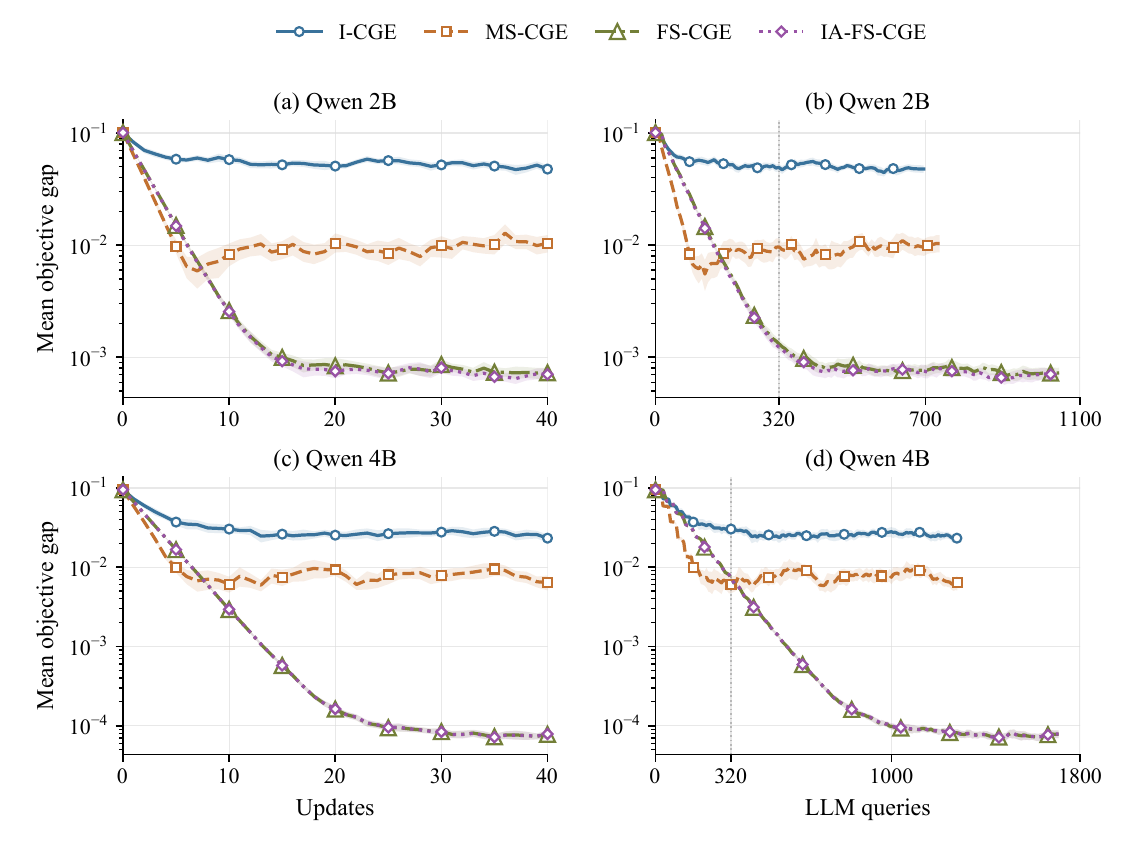}
    \caption{Software-onboarding optimization by updates (left) and LLM queries (right).}
    \label{fig:llm-onboarding}
\end{figure}

Figure~\ref{fig:llm-onboarding} reports the mean objective gaps during optimization for each method, with $c=0.025$, $n=4$, and $\theta_0=0.8$. The step sizes have been fine-tuned and selected as $(0.1,0.3,0.3,0.3)$ for I, MS, FS, and IA, respectively. Both subfigures use the same 40-update trajectories, plotted against updates and cumulative LLM queries. After 40 updates, FS and IA attain similar objective gaps, both substantially below those of I and MS. At a fixed budget of 320 calls, however, FS improves on MS under the 2B model, whereas the corresponding difference remains inconclusive under the 4B model.
This distinction arises because MS can complete more updates under small query budgets, while FS and IA require additional calls for residual continuations. As the budget increases, their greater gradient precision can offset this additional cost. At the fixed budgets specified in Appendix~\ref{app:llm-tuning}, IA achieves lower mean gaps than FS under the 2B model, but the adjusted paired comparisons remain inconclusive. Under the 4B model, the negligible difference between FS and IA may reflect similar loss impacts across events, which leave little scope for improvement through impact-aware allocation. Thus, IA may offer limited additional benefit over FS when impacts are comparable across events.

\section{Conclusion}\label{sec:conclusion}

LLM-embedded simulation optimization suffers from action noise that physical common random numbers cannot remove. Our score-access maximal coupling for finite-difference gradient estimators coordinates actions while preserving their marginals, giving an $O(c^{-1})$ variance bound under Lipschitz pre-switch stability. First-switch conditioning then averages over whether a switch occurs and yields uniformly bounded variance under pathwise stability. Its impact-aware variant allocates the first-switch event based on predicted event impacts to further improve the variance. These bounds quantify gradient accuracy and yield improved finite-sample bounds of projected stochastic approximation. The synthesized tests and LLM-embedded software-onboarding experiments support the variance bounds and their optimization performance. There are some open challenges for future directions.
For example, validation or alignment of LLM-based actions against human responses is needed to establish behavioral fidelity. Extensions to randomized directions require additional stability and error analysis.

\bibliographystyle{apalike}
\bibliography{reference}

\newpage

\begin{center}
    {\LARGE\bfseries Appendices}
\end{center}

\appendix
\numberwithin{theorem}{section}
\numberwithin{lemma}{section}
\numberwithin{proposition}{section}
\numberwithin{corollary}{section}
\numberwithin{assumption}{section}
\numberwithin{definition}{section}
\section{Supporting Results and Proofs}
\label{app:proofs}

When $\theta$ and $c$ are fixed in a proof, we may suppress these arguments in estimator notation.

\subsection{Independent-sampling variance floor}
\label{app:proof-indfloor}

\begin{theorem}
\label{thm:indfloor}
Fix $\theta\in\Theta$ and coordinate $j$. Let $\Delta L_{j,c}^{\mathrm{I}}=L_{j}^{\mathrm{I},+}-L_{j}^{\mathrm{I},-}$. For a behavioral event $e$, let $\mathcal H_{e,c}$ contain the shared physical primitives and every primitive behavioral seed except the one used to draw $A_e^+$. 
Assume there exist an event $e^*$, actions $a\ne b$, constants $c_1,p_0, q, \Delta_L>0$, and events $E_c\in\mathcal H_{e^*,c}$ such that $\inf_{0<c\le c_1}\Prob(E_c)\ge p_0$ and, on $E_c$,
\begin{enumerate}[label=(\roman*)]
\item $\Prob(A_{e^*}^+=a\mid \mathcal H_{e^*,c})\ge q$ and $\Prob(A_{e^*}^+=b\mid \mathcal H_{e^*,c})\ge q$.
\item with the randomness in $\mathcal H_{e^*,c}$ fixed, changing $A_{e^*}^+$ from $a$ to $b$ changes $\Delta L_{j,c}^{\mathrm{I}}$ by at least $\Delta_L$.
\end{enumerate}
Then there exists $v_0>0$, independent of $c$ and $n$, such that
\[
\Var\bigl(\widehat G_{j,n}^{\mathrm{I}}(\theta,c)\bigr)
\ge \frac{v_0}{nc^2}
\]
when $ 0<c\le c_1$.
\end{theorem}

\paragraph{Proof.}
Write $\mathcal H=\mathcal H_{e^*,c}$ and $Y=\Delta L_{j,c}^{\mathrm I}$ for simplicity. Conditional on $\mathcal H$, let $p_\ell$ be the probability of action $\ell$ and $y_\ell$ the resulting loss difference. All other primitive seeds are fixed, so
\[
\Var(Y\mid\mathcal H)
=\frac12\sum_{\ell,\ell'}p_\ell p_{\ell'}(y_\ell-y_{\ell'})^2
\ge p_a p_b(y_a-y_b)^2.
\]
On $E_c$, actions $a$ and $b$ each have probability at least $q$, and their corresponding loss differences are at least $\Delta_L$ apart. Hence $\Var(Y\mid\mathcal H)\ge q^2\Delta_L^2$ on $E_c$. The law of total variance and $\Prob(E_c)\ge p_0$ then give $\Var(Y)\ge p_0q^2\Delta_L^2$. Dividing by $4c^2$ and averaging $n$ independent replications yields
\[
\Var\!\left(\widehat G_{j,n}^{\mathrm I}(\theta,c)\right)
=\frac{\Var(\Delta L_{j,c}^{\mathrm I})}{4nc^2}
\ge\frac{p_0q^2\Delta_L^2}{4nc^2}.
\]
Hence $v_0=p_0q^2\Delta_L^2/4$ is uniform in $c$ and $n$.

\begin{corollary}
\label{cor:indmse}
Under Assumption~\ref{assump:smoothness} and the conditions of Theorem~\ref{thm:indfloor}, every sequence $c_n\downarrow0$ with $c_n\le c_1$ eventually satisfies
\[
\E\left[\big(\widehat G_{j,n}^{\mathrm{I}}(\theta,c_n)-\partial_j J(\theta)\big)^2\right]
\ge \frac{v_0}{n c_n^2}.
\]
In particular, with the canonical perturbation radius $c_n=\Theta ( n^{-1/4})$,
\[
\E\left[\big(\widehat G_{j,n}^{\mathrm{I}}(\theta,c_n)-\partial_j J(\theta)\big)^2\right]
=\Omega \left(n^{-1/2}\right).
\]
\end{corollary}

\paragraph{Proof.}
For every $c_n$, the bias--variance decomposition gives
\[
\E\bigl[(\widehat G_{j,n}^{\mathrm I}-\partial_jJ(\theta))^2\bigr]
\ge\Var(\widehat G_{j,n}^{\mathrm I})
\ge\frac{v_0}{nc_n^2}.
\]
Theorem~\ref{thm:indfloor} applies once $c_n\le c_1$. Substituting $c_n=\Theta(n^{-1/4})$ gives the stated $\Omega(n^{-1/2})$ lower bound.

\subsection{Preservation of the marginal laws}
\label{app:proof-pge-mean}

\begin{lemma}
\label{lemma:pge-mean-preservation}
Assume $\E[|L(\theta\pm ce_j,\boldsymbol\xi,\boldsymbol\omega)|]<\infty$ for the two standalone simulations. Under the exact-marginal construction in Definition~\ref{def:external-maximal}, the coupled trajectory losses for coordinate $j$ and replication $i$ satisfy
\[
\E[L_{j,i}^{\mathrm{MS},+}]=J(\theta+ce_j),
\qquad
\E[L_{j,i}^{\mathrm{MS},-}]=J(\theta-ce_j).
\]
Consequently,
\[
\E\!\left[\widehat G_{j,n}^{\mathrm{MS}}(\theta,c)\right] =D_{j,c}J(\theta) :=\frac{J(\theta+ce_j)-J(\theta-ce_j)}{2c}.
\]
\end{lemma}

\paragraph{Proof.}

Fix the physical primitives $\boldsymbol\xi$. Let $H_e^\pm$ be the paired pre-event histories, $C_e(h)=\Gamma_e(h,\zeta_e)$ the context, and $\Phi_e(h,a,\boldsymbol\xi)$ the history after action $a$ and the physical update. For bounded $f$, define the transition operator of the standalone $+$ simulation by
\[
(\mathcal T_ef)(h)
=\sum_{a\in\mathcal A}f\bigl(\Phi_e(h,a,\boldsymbol\xi)\bigr)P_{b,e}(a\mid C_e(h)).
\]
Exact marginality in \eqref{eq:maximal-properties} implies
\[
\E[f(H_{e+1}^+)\mid H_e^+,H_e^-,\boldsymbol\xi]=(\mathcal T_ef)(H_e^+).
\]
The right-hand side depends only on the $+$ history and gives its standalone transition at $\theta+ce_j$.
Taking conditional expectations given $(H_e^+,\boldsymbol\xi)$ leaves the same operator. The initial marginal law is unchanged, so induction preserves the full standalone history law at $\theta+ce_j$. Applying the integrable trajectory loss and averaging over $\boldsymbol\xi$ gives $\E[L_{j,i}^{\mathrm{MS},+}]=J(\theta+ce_j)$. The same argument for the $-$ side proves the central-difference identity.

\subsection{Proof of Theorem~\ref{thm:pge-main}}
\label{app:proof-pge-main}

Lemma~\ref{lemma:pge-mean-preservation} gives mean $D_{j,c}J(\theta)$. By Assumption~\ref{assump:smoothness},
\[
\left|D_{j,c}J(\theta)-\partial_jJ(\theta)\right|
\le\frac1{2c}\int_{-c}^{c}
\left|\partial_jJ(\theta+te_j)-\partial_jJ(\theta)\right|\,dt
\le\frac{L_Jc}{2}.
\]
For one pair, write $\Delta L_j=L_j^{\mathrm{MS},+}-L_j^{\mathrm{MS},-}$ and $N_j=N_{\mathrm{as},j}^{\mathrm{MS}}(\theta,c)$. Since $N_j^2\le mN_j$, Assumptions~\ref{assump:switch_dominance} and~\ref{assump:switch_sparsity} imply
\[
\E[(\Delta L_j)^2]
\le2L_f^2c^2+2C_{\mathrm{as}}^2\E[N_j^2]
\le2L_f^2c^2+2C_{\mathrm{as}}^2mC_{\mathrm{sp}}c^\kappa.
\]
Dividing this second-moment bound by $4c^2$ and averaging $n$ independent replications gives \eqref{eq:pge-component-variance}. Summing the coordinate bias--variance decompositions gives
\[
\E\|\widehat{\boldsymbol g}_n^{\mathrm{MS}}-\nabla J(\theta)\|^2
=\sum_{j=1}^d\left\{
\bigl(\E[\widehat G_{j,n}^{\mathrm{MS}}]-\partial_jJ(\theta)\bigr)^2
+\Var(\widehat G_{j,n}^{\mathrm{MS}})\right\},
\]
which proves \eqref{eq:pge-vector-mse}. The squared Euclidean norm has no cross-coordinate terms, so coordinate independence is unnecessary.

\begin{corollary}
\label{cor:pge-mse-rate}
Suppose Assumptions~\ref{assump:smoothness}, \ref{assump:switch_dominance}, and \ref{assump:switch_sparsity} hold.
Let $\bar\kappa:=\min\{\kappa,2\}$ and keep $d$ fixed. Choosing $c_n=\Theta\big(n^{-1/(4-\bar\kappa)}\big)$ gives
\[
\E\!\left[\left\|\widehat{\boldsymbol g}_{n}^{\mathrm{MS}}(\theta,c_n)-\nabla J(\theta)\right\|^2\right]
=O\!\left(d\,n^{-2/(4-\bar\kappa)}\right).
\]
\end{corollary}

\paragraph{Proof.}
With $d$ fixed, \eqref{eq:pge-vector-mse} bounds MSE by a constant multiple of
$d\{c^2+n^{-1}+n^{-1}c^{\kappa-2}\}$, with the constant independent of $n$ and $c$.
For $0<\kappa<2$, choose $c_n=a n^{-1/(4-\kappa)}$ for a fixed $a>0$.
Then $c_n^2=\Theta\bigl(n^{-2/(4-\kappa)}\bigr)$ and $n^{-1}c_n^{\kappa-2}=\Theta\bigl(n^{-2/(4-\kappa)}\bigr)$.
The remaining $n^{-1}$ term is no larger because $2/(4-\kappa)<1$.
For $\kappa\ge2$, choose $c_n=a n^{-1/2}$.
The squared-bias and smooth-path terms are $O(n^{-1})$, while
$n^{-1}c_n^{\kappa-2}=O(n^{-\kappa/2})=O(n^{-1})$.
Both choices eventually satisfy $c_n\le c_0$.
Combining the two cases gives the stated bound with $\bar\kappa=\min\{\kappa,2\}$.

\subsection{Proof of Proposition~\ref{prop:switchsparse}}
\label{app:proof-switchsparse}

The event $\{\sigma_j\ge e\}$ is $\mathcal F_{e-1}$-measurable. Conditional maximal coupling and Assumption~\ref{assump:preswitch} give
\[
\Prob(\sigma_j=e)
=\E[\ind\{\sigma_j\ge e\}\delta_{e,j}]
\le M_ec^\alpha.
\]
If the first switch occurs at $e$, at most $m-e+1$ events can contribute to the switch count. Thus
\[
N_{\mathrm{as},j}^{\mathrm{MS}}\le\sum_e(m-e+1)\ind\{\sigma_j=e\},
\]
and taking expectations proves the bound without restricting post-switch contexts.

Under the state-update conditions in \eqref{eq:preswitch-primitives}, set $K_1=\beta_0$ and $K_{e+1}=r_eK_e+\beta_e$. Induction gives $D_e\le K_ec$ over the finite horizon, even when some $r_e\ge1$. Combining the score bound with the softmax inequality below gives $\delta_{e,j}\le(L_sK_e^\alpha/\tau)c^\alpha$.

To verify the softmax inequality, let $\sigma_\tau=\operatorname{softmax}_\tau$ and $z(t)=y+t(x-y)$. Then
\[
\frac{d}{dt}\sigma_{\tau,a}(z(t))
=\frac{\sigma_{\tau,a}(z(t))}{\tau}
\left\{(x_a-y_a)-\sum_b\sigma_{\tau,b}(z(t))(x_b-y_b)\right\}.
\]
Both the coordinate difference and its probability-weighted mean are bounded by $\|x-y\|_\infty$. Summing absolute derivatives and integrating yields
\[
\TV(\sigma_\tau(x),\sigma_\tau(y))\le\|x-y\|_\infty/\tau.
\]

\subsection{Proof of Proposition~\ref{prop:first-order-switch}}
\label{app:proof-first-order}

Keep $\boldsymbol\xi$ fixed while varying $\theta$. Differentiability and the finite action space imply, almost surely,
\[
\frac{\TV(p_1^+,p_1^-)}{c}
=\sum_{a\in\mathcal A}\left|\frac{p_1^+(a)-p_1^-(a)}{2c}\right|
\longrightarrow\sum_{a\in\mathcal A}
|\partial_jP_{b,1}(a\mid C_1^\theta)|.
\]
Here $p_1^\pm(a)=P_{b,1}(a\mid C_1^{\theta\pm ce_j})$. Maximal coupling makes the first-event switch probability equal to this TV distance, and the total count is at least that switch indicator. Fatou's lemma therefore gives
\[
\liminf_{c\downarrow0}\frac{\E[N_{\mathrm{as},j}^{\mathrm{MS}}(\theta,c)]}{c}
\ge\E\left[\sum_{a\in\mathcal A}
|\partial_jP_{b,1}(a\mid C_1^\theta)|\right]>0.
\]
The sensitivity condition makes the expectation positive, whether finite or infinite. Hence the count is $\Omega(c)$, ruling out an $O(c^\kappa)$ bound with $\kappa>1$.

\paragraph{Variance of a paired finite difference.}
For a one-event binary model with loss $zA$, $z\ne0$, let $p^\pm=\Prob(A^\pm=1)$ and $\Delta p=p^+-p^-$. Every exact-marginal coupling satisfies
\[
\Var\!\left(\frac{z(A^+-A^-)}{2c}\right)
\ge\frac{z^2}{4c^2}\left(|\Delta p|-(\Delta p)^2\right),
\]
because $\E[(A^+-A^-)^2]\ge|\Delta p|$. Maximal coupling attains this bound, which is $\Omega(c^{-1})$ when $|\Delta p|=\Theta(c)$. First-switch conditioning averages over the no-switch and first-switch cases and can avoid this variance growth.

\subsection{Proof of Theorem~\ref{thm:fs-conditioning}}
\label{app:proof-fs-conditioning}

Fix $\theta,j,c$ and let $\mathcal B$ record the physical primitives and sampled common-action backbone. Conditional on the physical primitives, a reached prefix $a_1,\ldots,a_{e-1}$ has backbone probability $\prod_{r<e}h_r(a_r)/\rho_r$. Multiplying by its first-switch weight gives
\[
\left(\prod_{r<e}\frac{h_r(a_r)}{\rho_r}\right)w_e
=\left(\prod_{r<e}h_r(a_r)\right)\delta_e,
\]
which is the MS-CGE probability of this prefix followed by a first switch at $e$. The residual draw and subsequent continuation also follow the MS-CGE law. Similarly, multiplying a complete backbone's probability by $W_m$ gives its no-switch probability $\prod_rh_r(a_r)$. If an overlap is zero, the backbone stops before the common action is drawn and $W_m=0$.

Let $\Delta^{\mathrm{MS}}$ be the MS loss difference and $\Delta_e$ a residual difference from the saved prefix at $e$. The no-switch and first-switch cases partition the sample space, so for bounded $f$,
\[
\E[f(\Delta^{\mathrm{MS}})]
=\E\left[W_mf(\Delta_0)+\sum_ew_e\E[f(\Delta_e)\mid\mathcal B]\right].
\]
Fresh residual draws are independent of future backbone draws. For $S>0$, selecting $R$ with probabilities $w_e/S$ expresses the sum as $S\E[f(\Delta_R)\mid\mathcal B]$. The loss bound allows us to take $f(x)=x$ or $f(x)=x^2$. The first choice proves $\E[G_j^{\mathrm{FS}}]=D_{j,c}J(\theta)$.

For the variance identity, set $\Delta_R=\Delta_0$ when $S=0$ and, conditional on $(\mathcal B,R,\Delta_R)$, draw $I_S\sim\mathrm{Bernoulli}(S)$. Then $Y=(1-I_S)\Delta_0+I_S\Delta_R$ has the law of $\Delta^{\mathrm{MS}}$, with
\[
\E[Y\mid\mathcal B,R,\Delta_R]=2cG_j^{\mathrm{FS}},\qquad
\Var(Y\mid\mathcal B,R,\Delta_R)=S(1-S)(\Delta_R-\Delta_0)^2.
\]
The law of total variance, divided by $4c^2$, proves \eqref{eq:fs-variance-identity}. When $S=0$, the branch is omitted and its variance contribution is zero.

\subsection{Proof of Theorem~\ref{thm:fs-rate}}
\label{app:proof-fs-rate}

We first show that the pre-switch bound assumed for MS-CGE also holds on the common-action backbone. Let $\nu_e$ be the joint subprobability measure of physical inputs and backbone prefixes that reach event $e$, and let $\widetilde\nu_e$ be the corresponding MS-CGE measure restricted to paths with no earlier switch. Conditional on the physical inputs, their prefix masses are $\prod_{r<e}h_r(a_r)/\rho_r$ and $\prod_{r<e}h_r(a_r)$, respectively. Thus $d\widetilde\nu_e=W_{e-1}\,d\nu_e$. Since $W_{e-1}>0$ on reached prefixes, both laws assign zero mass to the same sets there. Thus Assumption~\ref{assump:preswitch} with $\alpha=1$ holds almost surely on the backbone. Summing over reached events gives
\[
S=\sum_eW_{e-1}\delta_e\le\sum_e\delta_e\le Mc.
\]
If an overlap is zero, its first-switch weight is included before the backbone stops, and all later weights are zero.
All backbone actions agree, while a residual continuation can contain at most $m$ switches. Assumption~\ref{assump:switch_dominance} gives $|\Delta_0|\le L_fc$ and $|\Delta_R|\le B=L_fc_0+C_{\mathrm{as}}m$, so
\[
|G_j^{\mathrm{FS}}|
\le\frac{W_mL_fc+SB}{2c}
\le\frac{L_f+MB}{2}=C_{\mathrm{FS}}.
\]
The bound also holds when an overlap is zero because $W_m\Delta_0=0$. Averaging $n$ independent replications gives variance at most $C_{\mathrm{FS}}^2/n$. Combining this with the bias bound in Appendix~\ref{app:proof-pge-main} finishes the proof of \eqref{eq:fs-mse}. 

\subsection{Proof of Proposition~\ref{prop:ia-fs-allocation}}
\label{app:proof-ia-allocation}

Fix the predictor and condition on $\mathcal B$. Write $\mu_e=\E[\Delta_e\mid\mathcal B]$, $\mathcal M_e=\E[\Delta_e^2\mid\mathcal B]$, and use $T=\sum_e w_e\mu_e$ locally in this proof. Each branch reuses its saved prefix and physical inputs with fresh residual draws, so its law is independent of $q$. Averaging over branch selection gives
\[
\E[2cG_j^{\mathrm{IA}}(\theta,c;q)\mid\mathcal B]
=W_m\Delta_0+\sum_{e:q_e>0}q_e\frac{w_e}{q_e}\mu_e
=W_m\Delta_0+T.
\]
Cancellation of $q_e$ makes this mean independent of the selection probabilities. Appendix~\ref{app:proof-fs-conditioning} then gives the central-difference mean. The bounds $|\Delta_e|\le B$ and $\sum_e w_e=S\le1$ ensure integrability. When $S=0$, only the backbone term remains, and we define the integrands in \eqref{eq:ia-fs-oracle-gain} and \eqref{eq:ia-fs-allocation-penalty} to be zero.

When the estimator has finite variance, conditioning on $\mathcal B$ fixes the no-switch term and gives
\begin{equation}
\Var(G_j^{\mathrm{IA}}(\theta,c;q)\mid\mathcal B)
=\frac{1}{4c^2}\left(\sum_{e:q_e>0}\frac{w_e^2\mathcal M_e}{q_e}-T^2\right).
\label{eq:ia-fs-conditional-variance}
\end{equation}
Only the sum depends on $q$. By Cauchy--Schwarz, its minimum is $Z^2=(\sum_e w_e\sqrt{\mathcal M_e})^2$, attained by $q^\star$ when $Z>0$ and zero-contribution branches may be omitted. Such an omission is valid because $\mathcal M_e=0$ implies $\Delta_e=0$ conditionally almost surely. When $Z>0$, requiring every positive-weight branch to receive positive probability makes $Z^2$ an unattained infimum if any such branch has $\mathcal M_e=0$. It is approached as the probabilities of all such zero-moment branches tend to zero. The clipped implementation gives positive probability to every positive-weight branch. When $Z=0<S$, all residual contributions vanish. We then set $q^\star=q^0$ and $Z^2\chi^2=0$, so $q^0$ also has zero residual variance.

For FS-CGE, the sum in \eqref{eq:ia-fs-conditional-variance} is $S\sum_e w_e\mathcal M_e$. Subtracting the oracle value $Z^2$ from the FS and general-allocation sums gives
\[
S\sum_e w_e\mathcal M_e-Z^2
=S^2\Var_{e\sim q^0}(\sqrt{\mathcal M_e}),\qquad
\sum_{e:q_e>0}\frac{w_e^2\mathcal M_e}{q_e}-Z^2
=Z^2\chi^2(q^\star\Vert q).
\]
The variance of the common conditional mean cancels in the law of total variance. Taking expectations of these identities and dividing by $4c^2$ proves \eqref{eq:ia-fs-oracle-gain} and \eqref{eq:ia-fs-allocation-penalty}. The branch bounds ensure finite second moments for FS-CGE and the oracle at fixed $c>0$.

For the single-checkpoint implementation, conditioning on $\mathcal B$ and $H$ fixes the reservoir weights $v_e$. Event $e$ becomes the final checkpoint if it enters and is never replaced. For $v_e>0$ and $V_e=\sum_{r\le e}v_r$, this probability is
\[
\frac{v_e}{V_e}\prod_{r>e:v_r>0}\left(1-\frac{v_r}{V_r}\right)
=\frac{v_e}{V_e}\prod_{r>e:v_r>0}\frac{V_{r-1}}{V_r}
=\frac{v_e}{V_m}.
\]
Averaging over $H$ gives the selection probabilities $q^\gamma$ in \eqref{eq:ia-fs-proposal}. Using these probabilities in the correction gives \eqref{eq:ia-fs-mixture-weight}. Constant multipliers recover FS-CGE.

\subsection{Proof of Corollary~\ref{cor:ia-fs-rate}}
\label{app:proof-ia-rate}

Appendix~\ref{app:proof-fs-rate} gives $S\le Mc$ on the backbone.

The loss assumption gives $|\Delta_0|\le L_fc$ and $|\Delta_R|\le B$. Clipping $u_e$ gives $q_e^\gamma\ge q_e^0/\Lambda$, which bounds the residual importance factor by $w_R/q_R^\gamma\le\Lambda S$. Combining these bounds yields
\[
|G_j^{\mathrm{IA}}(\theta,c;q^\gamma)|
\le\frac{W_mL_fc+\Lambda SB}{2c}
\le\frac{L_f+\Lambda MB}{2}=C_{\mathrm{IA}}.
\]
For a fixed predictor, this bound gives coordinate variance at most $C_{\mathrm{IA}}^2/n$ after averaging $n$ independent replications. Averaging over training randomness preserves the bound because the conditional mean $D_{j,c}J(\theta)$ is unchanged. Adding squared bias and summing over coordinates proves \eqref{eq:ia-fs-mse}. Taking $c_n=\Theta(n^{-1/2})$ gives $O(d/n)$ MSE.

\subsection{Proof of Theorem~\ref{thm:sample_complexity}}
\label{app:proof-of-cge-complexity}

Let $\mathcal F_k$ contain the history before update $k$, including the fixed IA predictor, and write $\E_k=\E[\cdot\mid\mathcal F_k]$. Fresh simulation batches preserve the central-difference mean. For each method, the bias argument in Appendix~\ref{app:proof-pge-main} gives a squared conditional bias bound of $b(c)^2=dL_J^2c^2/4$. For I-CGE, $N_{\mathrm{as},j}^{\mathrm I}\le m$ yields
\[
\E_k\|\widehat{\boldsymbol g}_k^{\mathrm I}-\E_k\widehat{\boldsymbol g}_k^{\mathrm I}\|^2
\le\frac{dL_f^2}{2n}+\frac{dC_{\mathrm{as}}^2m^2}{2nc^2}.
\]
For MS-CGE, $N^2\le mN$ and Assumption~\ref{assump:switch_sparsity} replace $m^2$ by $mC_{\mathrm{sp}}c^\kappa$. The FS and IA bounds are $dC_{\mathrm{FS}}^2/n$ and $dC_{\mathrm{IA}}^2/n$, respectively. Let $v(c,n)$ denote the resulting bound on the conditional variance of the gradient estimate.

Set $e_k=\widehat{\boldsymbol g}_k-\nabla J(\theta_k)$ and $\widehat{\mathcal G}_{\eta,k}=(\theta_k-\theta_{k+1})/\eta$. Bias--variance decomposition then gives $\E_k\|e_k\|^2\le b(c)^2+v(c,n)$. Projection optimality gives
\[
\langle\widehat{\boldsymbol g}_k,\widehat{\mathcal G}_{\eta,k}\rangle\ge\|\widehat{\mathcal G}_{\eta,k}\|^2.
\]
Smoothness, Young's inequality and $\eta\le1/(2L_J)$ imply
\[
\begin{aligned}
J(\theta_{k+1})
&\le J(\theta_k)-\eta\langle\widehat{\boldsymbol g}_k,\widehat{\mathcal G}_{\eta,k}\rangle
+\eta\langle e_k,\widehat{\mathcal G}_{\eta,k}\rangle
+\frac{L_J\eta^2}{2}\|\widehat{\mathcal G}_{\eta,k}\|^2\\
&\le J(\theta_k)-\frac\eta4\|\widehat{\mathcal G}_{\eta,k}\|^2+\frac\eta2\|e_k\|^2.
\end{aligned}
\]
Nonexpansiveness of projection gives
$\|\widehat{\mathcal G}_{\eta,k}-\mathcal G_\eta(\theta_k)\|\le\|e_k\|$.
Consequently,
\[
J(\theta_{k+1})\le J(\theta_k)-\eta\|\mathcal G_\eta(\theta_k)\|^2/8+3\eta\|e_k\|^2/4.
\]
Taking expectations, summing over $k<K$, and using $J(\theta_K)\ge J_{\inf}$ yields
\begin{equation}
\min_{0\le k<K}\E\|\mathcal G_\eta(\theta_k)\|^2
\le\frac{8(J(\theta_0)-J_{\inf})}{\eta K}+6b(c)^2+6v(c,n).
\label{eq:generic-projected-bound}
\end{equation}

With fixed $\eta$, choose $K=\Theta(\epsilon^{-1})$ and $c=\Theta(\sqrt{\epsilon/d})$, so the first two terms in \eqref{eq:generic-projected-bound} are $O(\epsilon)$. To obtain the same order for the variance term, take $n=\Theta(d^2\epsilon^{-2})$ for I-CGE. For MS-CGE, take $n=\Theta(d^{2-\kappa/2}\epsilon^{-2+\kappa/2})$ when $0<\kappa<2$ and $n=\Theta(d\epsilon^{-1})$ when $\kappa\ge2$. FS-CGE and IA-FS-CGE also require $n=\Theta(d\epsilon^{-1})$. Choosing the constants appropriately makes each term at most $\epsilon/3$, with $c\le\min\{c_0,1\}$. Multiplying the batch size by $dK$ gives the coordinate-replication bounds.

\section{Numerical Settings and Parameter Selection}
\label{app:numerical-settings}

\subsection{Evaluation and training protocols}
\label{app:numerical-measures}

We enumerate the finite simulation tree to compute exact objectives, gradient bias, variance and MSE. Optimization gaps $J(\theta)-J(\theta^\star)$ are evaluated at each run's last completed iterate and then averaged across runs, where $\theta^\star$ minimizes $J$. The count $N_{\mathrm{as}}^{\mathrm{branch}}$ includes the forced first switch and later action disagreements. It is zero when no residual continuation is generated.

One evaluation returns a behavioral probability vector at a given context and perturbed policy. The LLM experiments use cached scores, counting each probability-vector evaluation as a query. Online budgets count backbone and residual evaluations but exclude predictor training and exact diagnostics. Before each update, I and MS reserve $2mn$ calls, while FS and IA reserve $(4m-2)n$ by reusing the selected event's probability vectors. Only actual calls are charged. Fixed-budget runs have no update limit and report the last completed iterate. Physical inputs are indexed by event and update to keep methods paired across different update counts. Predictor training, backbone sampling, checkpoint selection and residual sampling use separate random streams.

At each radius, we train independent predictors and fix them before evaluation. Each of the $P$ training attempts samples a policy uniformly from 17 equally spaced values in the policy interval and a target event uniformly from the four events. If the event is reached and permits a first switch, a forced continuation supplies the label $\Delta_e^2$. Otherwise, no label is obtained, but the calls are counted. Inputs are policy, event, profile and pre-action state. We average training labels at each policy value, linearly interpolate these averages and take the square root to obtain the common RMS scale. Dividing the predicted RMS by this scale and clipping to $[0.5,2]$ gives $u_e$. If fitting fails, $u_e=1$.

Each model uses twenty independent predictor fits per radius, reused across figures at shared radii. Optimization uses 32 paired runs per fit and method, giving 640 runs per method. Bootstrap resampling draws fits as blocks, retaining each fit's 32 paired runs together. Optimization bands are pointwise $95\%$ intervals based on 4,000 resamples. Fixed-budget IA--FS comparisons in the LLM models use 20,000 resamples and a Bonferroni correction for twelve comparisons, giving nominal $95\%$ simultaneous coverage. The separate numerical supplement documents the full comparison set and training costs.

\subsection{Synthesized model and settings}
\label{app:controlled-settings}
\label{app:controlled-tuning}

At each of four interactions, $A_e=1$ denotes satisfaction and $A_e=0$ dissatisfaction. We set $A_0=0$ and draw a persistent customer type $Z\in\{-0.4,0.4\}$ with equal probabilities. For $\theta\in[0.1,1.9]$, the satisfaction probability and loss are
\[
p_\theta(A_{e-1},Z)=\ell\bigl(1.5\theta-1+0.8(2A_{e-1}-1)+Z\bigr),\qquad
L(\theta)=\frac{\theta^2}{2}+\sum_{e=1}^4g_e(1-A_e),
\]
where $\ell(x)=(1+\exp(-x))^{-1}$ and $(g_1,g_2,g_3,g_4)=(0.8,0.4,0.4,0.4)$. Thus dissatisfaction at the first interaction costs twice as much as at any later interaction.

The kernel is defined on the policy interval's $0.2$-neighborhood for two-sided perturbations. Before the first switch, customer types and action histories agree, so $\ell'\le1/4$ gives $\delta_e\le0.75c$. Thus $M_e=0.75$, $\alpha=1$ and $\E[N_{\mathrm{as}}^{\mathrm{MS}}]\le7.5c$. The loss bound $|\Delta L|\le3.8c+0.8N_{\mathrm{as}}$ gives $L_f=3.8$, $C_{\mathrm{as}}=0.8$, $c_0=0.2$, $M=3$, $B=3.96$ and $C_{\mathrm{FS}}=7.84$. Clipping gives $\Lambda=4$ and $C_{\mathrm{IA}}=25.66$. The first-event satisfaction probability has a positive policy derivative, as required by Proposition~\ref{prop:first-order-switch}.

For a fixed customer type $z$, let $p_a=p_\theta(a,z)$ and $v_e=\Prob(A_e=1\mid Z=z)$. The recursion $v_e=p_0+(p_1-p_0)v_{e-1}$, starting from $v_0=0$, gives
\[
J(\theta)=\frac{\theta^2}{2}+2-\frac{1}{2}\sum_z\sum_{e=1}^4g_ev_e(\theta,z).
\]
Differentiating gives $J'$ and $J''$. The optimum is $\theta^\star\approx0.808829$, with $J(\theta^\star)\approx1.400605$. The bounds $|p_a'|\le1.5/4$ and $|p_a''|\le1.5^2/(6\sqrt3)$ imply
\[
|v_e''|\le1.5^2e/(6\sqrt3)+1.5^2e(e-1)/16
\]
and $|J''|<3.077628$, so we use $L_J=3.08$.

The variance and switch comparisons use $\theta=0.2$ and radii $c=0.2\times2^{-i}$, $i=0,\ldots,6$. For the variance--evaluation-cost comparison, we use 17 equally spaced policies in $[0.1,1.9]$ at $c=0.05$. We average the IA/FS product ratios over these policies and then over fits, and report one minus this average. The MSE experiment uses $n=16\times2^i$, $i=0,\ldots,8$, with radii $0.4n^{-1/4}$, $0.4n^{-1/3}$ and $0.4n^{-1/2}$ for I, MS and FS, respectively. IA uses the FS schedule. All four curves use exact moments, with IA averaged over twenty fits. Optimization tests use $\theta_0=0.2$, $c=0.05$, $n=8$, $\eta=0.1$ and 80 updates for each method.

At each training policy, the IA predictor averages squared loss differences by event, customer type, and pre-action state. To stabilize estimates from small groups, we combine event means with the overall mean and group means with their event means, giving the broader mean a weight equivalent to 32 labels at both levels. Predictions are linearly interpolated between training policies. We use $P=32{,}768$ attempts and $\gamma=1$, which requires $262{,}144$ training calls per fit.

\subsection{Software-onboarding model and settings}
\label{app:llm-settings}
\label{app:llm-tuning}
\label{app:llm-validation}

A user receives up to four invitations to set up a software account and has one of two equally likely profiles: a time-constrained freelance business owner or a community-group coordinator wary of unfamiliar settings. The profile remains fixed throughout onboarding. State $s$ records the profile, invitation round and previous deferral or technical failure. Basic invitations offer self-service setup, while assisted invitations add staff guidance. Prices, permissions and the twenty-minute setup task are the same under both invitation types. Users choose whether to attempt setup or defer. Each attempt succeeds with probability $0.8$, and successful activation ends the user's onboarding process.

We run Qwen3.5-2B and Qwen3-4B-Instruct-2507 locally in MLX with 4-bit weights and fixed prompts. For each state and invitation type, we obtain the evaluation-mode next-token logits for A and B under both label assignments: attempt=A and defer=B, followed by the reverse assignment. We align the logits by action and average the two attempt-minus-defer differences. Applying the logistic function with temperature 1 for 2B and 5 for 4B gives the attempt probabilities $q_{\mathrm{basic}}(s)$ and $q_{\mathrm{assist}}(s)$. Each profile has one initial state and two pre-action states per later round, giving fourteen active states in total. All scores are cached before optimization.

The assistance probability $\theta$ gives $p_\theta(s)=(1-\theta)q_{\mathrm{basic}}(s)+\theta q_{\mathrm{assist}}(s)$ because invitation type affects only the current action. Every estimator uses $(p_\theta,1-p_\theta)$ and the loss
\[
L(\theta)=\frac1{30}\left(8\theta^2+N_{\mathrm{contact}}+
30\ind\{\text{not activated after round }4\}\right),
\]
where $N_{\mathrm{contact}}$ is the number of invitations sent and the capacity cost is charged upfront. Shared physical success draws keep the paired states identical before the first action switch. Hence $\delta_e=2c|q_{\mathrm{assist}}(s_e)-q_{\mathrm{basic}}(s_e)|$ satisfies pre-switch stability with $M_e=2\max_s|q_{\mathrm{assist}}(s)-q_{\mathrm{basic}}(s)|$ and $\alpha=1$. Capacity costs differ by at most $0.96c$. The combined contact and failure costs differ by at most 1.1 and are equal when no switch occurs, giving $|\Delta L|\le0.96c+1.1N_{\mathrm{as}}$.

Policy centers lie in $[0.1,0.9]$. Dynamic programming gives quartic objectives and exact derivatives, with $(\theta^\star,J(\theta^\star))\approx(0.215,0.217)$ for 2B and $(0.204,1.109)$ for 4B. The gradient-moment table uses $\theta=0.8$ and $c=0.025$, while the radius plots use $\theta=0.8$ and $c=0.05\times2^{-i}$, $i=0,\ldots,4$. The central-difference bias is exactly $c^2J'''(\theta)/6$.

Optimization starts at $\theta_0=0.8$ with $c=0.025$ and $n=4$. We select step sizes from $\eta\in\{0.01,0.03,0.1,0.3\}$ using 12 development runs per candidate and the final objective within a budget of 320 evaluations. We set $L_J$ 1\% above the maximum $|J''|$ on $[0.05,0.95]$. Every candidate satisfies $\eta\le1/(2L_J)$. The selected step sizes are $(0.1,0.3,0.3)$ for I, MS and FS in both models, and IA uses the FS value. The iteration comparison uses 40 updates. Separate fixed-budget runs use 320, 1,280 and 5,120 evaluations. At 320 evaluations, the FS-minus-MS and MS-minus-I comparisons use 4,000 paired block-bootstrap resamples and $97.5\%$ marginal intervals to obtain nominal $95\%$ joint coverage within each model.

IA uses ridge regression with penalty 8 and an unpenalized intercept in both LLM models. Regressors include indicators for events after the first, the policy rescaled to $[-1,1]$, its square, and their event interactions. A profile code $\pm1$ and an indicator of prior technical failure enter directly and through interactions with the rescaled policy. Labels are divided by their overall training mean before fitting. Normalized predictions are bounded below by $0.001$ and rescaled before computing RMS multipliers. We use $(P,\gamma)=(32{,}768,1)$ for 2B and $(8{,}192,0.25)$ for 4B.

\end{document}